\documentclass[a4paper,12pt]{article}
\usepackage{ucs}
\usepackage[utf8x]{inputenc}
\usepackage{amsfonts} 
\usepackage{amsmath} 
\usepackage[dvips]{graphicx}

\usepackage{caption} 
\usepackage{subcaption} 
\usepackage{enumerate}
\usepackage{setspace}
\usepackage[toc,page]{appendix}
\usepackage{soul}

\usepackage{color}
\title{}
\author{}
\date{}

\newtheorem{theorem}{Theorem}[section]

\newtheorem{proposition}[theorem]{Proposition}
\newtheorem{corollary}[theorem]{Corollary}

\newenvironment{@abssec}[1]{%
\if@twocolumn
\section*{#1}%
\else
\vspace{.05in}\footnotesize
\parindent .2in
{\upshape\bfseries #1. }\ignorespaces
\fi}
{\if@twocolumn\else\par\vspace{.1in}\fi}

\newcommand{\qed}{\nobreak \ifvmode \relax \else
\ifdim\lastskip<1.5em \hskip-\lastskip
\hskip1.5em plus0em minus0.5em \fi \nobreak
\vrule height0.75em width0.5em depth0.25em\fi}

\newcommand{\Q}{\mbox{WB}}

\newcommand{\torus}{{\mathbb{T}^d}}
\newcommand{\TTT}{{\mathbb{T}}}
\newcommand{\Torus}{{\mathbb{T}^2}}

\newcommand{\R}{{\mathbb{R}}}

\newcommand{\mR}{\mathbb{R}}

\newcommand{\Path}{}
\newcommand{\Xx}{\ref{eqn:plusrho}}

\newcommand\keywordsname{Key words}
\newcommand\AMSname{AMS subject classifications}

\begin{document}
\title{Solving the Babylonian Problem of quasiperiodic rotation rates} 

\author{Suddhasattwa Das\footnotemark[1], \and Yoshitaka Saiki\footnotemark[2]\ \footnotemark[3]\ \footnotemark[5],
\and Evelyn Sander\footnotemark[4], \and James A Yorke\footnotemark[5]}
\footnotetext[1]{current address: Courant Institute of Mathematical Sciences, New York University}
\footnotetext[2]{Graduate School of Commerce and Management, Hitotsubashi University}
\footnotetext[3]{JST, PRESTO}
\footnotetext[4]{Department of Mathematical Sciences, George Mason University}
\footnotetext[5]{University of Maryland, College Park}

\footnotetext[1]{Department of Mathematics, University of Maryland, College Park}

\date{\today}
\maketitle

\renewcommand*\contentsname{Contents}

\begin{abstract} A trajectory $\theta_n := F^n(\theta_0), n = 0,1,2, \dots $ 
is quasiperiodic if the trajectory lies on and is dense in 
some $d$-dimensional torus, and there is a choice of coordinates on the torus $\torus$ 
for which $F$ has the form $F(\theta) = \theta + \rho\bmod1$
for all $\theta\in\torus$ and for some $\rho\in\torus$. 
(For $d>1$ we always interpret $\bmod1$ as being applied to each coordinate.)
There is an ancient literature on computing three rotation rates for the Moon.
However, for $d>1$, the choice of coordinates that yields the form $F(\theta) = \theta + \rho\bmod1$ 
is far from unique and the different choices yield a huge choice of coordinatizations $(\rho_1,\cdots,\rho_d)$ of $\rho$, 
and these coordinations are dense in $\torus$. Therefore instead one defines the rotation rate $\rho_\phi$ from the perspective of a 
map $\phi:T^d\to S^1$. This is in effect the approach taken by the Babylonians and we refer to this approach as the ``Babylonian Problem''.
However, even in the case $d=1$ there has been no general
method for computing $\rho_\phi$ given only the sequence $\phi(\theta_n)$, 

though there is a literature dealing with special cases. Here we
present our {\em Embedding continuation method} for computing
  $\rho_\phi$ from the image $\phi(\theta_n)$ of a trajectory. 
  
It is based on the Takens
Embedding Theorem and the Birkhoff Ergodic Theorem. 
\end{abstract}

{\bf Keywords:} Quasiperiodic, Birkhoff Ergodic Theorem, Rotation number, Rotation rate, Takens Embedding Theorem,
Circular Planar Restricted 3-Body Problem, CR3BP

\section{Introduction}\label{sec:introduction}

The goal of this paper is to show how to compute a rotation rate of a quasiperiodic
discrete-time trajectory. We begin with a motivating historical example, followed by 
a broad overview of our approach to determining rotation rates.

Rotation rates and quasiperiodicity have been studied for millennia; namely, 
the Moon's orbit has three periods whose approximate values
were found 2500 years ago by the Babylonians~\cite{Goldstein}. Although
computation of the periods of the Moon is an easy problem today, we use
it to give context to the problems we investigate. The Babylonians found that the 
periods of the
Moon - measured relative to the distant stars - are approximately 27.3
days (the sidereal month), 8.85 years for the rotation of the apogee
(the local maximum distance from the Earth), and 18.6 years for the
rotation of the intersection of the Earth-Sun plane with the Moon-Earth
plane. They also measured the variation in the speed of the Moon through
the field of stars, and the speed is inversely correlated with the
distance of the Moon. They used their results to predict
eclipses of the Moon, which occur only when the Sun, Earth and Moon are
sufficiently aligned to allow the Moon to pass through the shadow of the
Earth. How they obtained their estimates is not fully understood but it was through
observations of the trajectory of the Moon through the distant stars in the sky.
In essence they viewed the Moon projected onto the two-dimensional space
of distant stars.   We too  work with quasiperiodic motions which
have been projected into one or two dimensions. 

The Moon has three periods because the Moon's orbit
is basically three-dimensionally quasiperiodic, traveling on a
three-dimensional torus $\TTT^3$ that is embedded in six
(position+velocity) dimensions. The torus is topologically the product
of three circles, and the Moon has an (average) rotation {\em rate} 
--  i.e. the reciprocal of 
the period -- along each
of these circles. While the Moon's orbit has many intricacies, one can
capture some of the subtleties by approximating the Sun-Earth-Moon
system as three point masses using Newtonian gravitational laws. This leads to the
study of the Moon's orbit as a circular restricted three-body problem
(CR3BP) in which the Earth travels on a circle about the Sun and the
Moon has negligible mass. Using rotating coordinates in which the Earth
and Sun are fixed while the Moon moves in three-dimensions, the orbit
can thus be approximated by the above mentioned three-dimensional torus
$\TTT^3$ in $\mathbb{R}^6$. Such a model ignores several factors including long-term tidal forces
and the small influence of the other planets.

As another  motivating example,  the direction $\phi$ of Mars from the Earth (viewed against the backdrop of the fixed stars)
does not change monotonically. This apparent non-monotonic movement is
called ``retrograde motion.'' 
Now imagine that exactly once each year the
direction $\phi$ is determined. How do we determine the rotation rate
of Mars compared with an Earth year from such data? 

This kind of problem has
been unsolved in full generality even for images of one-dimensional quasiperiodic maps. 

This paper considers a setting more general than just  the Moon or on Mars, although both
give good illustrations of our setting. 

For typical
discrete-time dynamical systems, it is conjectured that the
three kinds of recurrent motions that are likely to be seen in a dynamical 
system are periodic
orbits, chaotic orbits, and quasiperiodic orbits
~\cite{sander:yorke:15}. Starting with a
$d$-dimensional quasiperiodic orbit on a torus $\torus$ for some $d$ and
a map $\phi:\torus\to S^1$, we  establish a new method for computing
rotation rates from a discrete-time quasiperiodic orbit. By discrete time, we mean that the
trajectory observations are a discrete sequence $\phi_n, n =
0,1,2,\cdots$, as for example when a Poincar\'e return map is used for
the planar circular restricted three-body problem (CR3BP), and $\phi_n$
is 
the angle of the image of the trajectory as measured from the perspective of some reference point 

at the $n^{th}$ time the trajectory crosses some specified
Poincar\'e surface.

In the rest of this introduction, we give a non-technical 
summary of our results, ending with a comparison to previous work on this topic.
We then proceed with a more technical parts of the paper, in which we describe our 
methods and results in detail and give numerical examples for which we compute rotation rates. 

{\bf Quasiperiodicity defined.} 
Let $\torus$ be a $d$-dimensional torus. 
A quasiperiodic orbit is an orbit that is dense on 
a $d$-dimensional torus and such that there  exists a choice of coordinates 
$\theta \in\torus
:= [0,1]^d \bmod1$ (where $\bmod1$ is applied to each coordinate) for
the torus such that the dynamics on the orbit are given by the map
\begin{equation}\label{eqn:plusrho}
\theta_{n+1} := F(\theta_n) = \theta_n + \rho\bmod1
\end{equation}

for some 
{\bf rotation vector} 
$\rho\in\torus$ where the coordinates $\rho_i$ of the $\rho$ are irrational
and {\bf rationally independent}, 
i.e. if $a_k$ are rational numbers for $k=1,\cdots,d$ for which 
$a_1\rho_1+\cdots+a_d\rho_d = 0$, then $a_k=0$ for all $k=1,\cdots,d$. 
We will say such a 
rotation vector $\rho$ is {\bf irrational}. 

{\bf The Babylonian Problem.} 
One might imagine that our goal would be to compute $\rho$ in Eq. \Xx\ from whatever knowledge we could obtain about the torus $\torus$. Although the Babylonians did not know about three-dimensional tori, they none the less obtained three meaningful rotation rates. 
To abstract their situation, 
we assume there is a smooth map $\psi:\torus\to M$ where $M$ is a manifold, usually of dimension $1$ or $2$. The {\bf Babylonian Problem} 
is to compute a rotation rate $\rho_\psi$ from knowledge of the projection of a trajectory. We assume we only have the values $\psi_n$ of $\psi$ at a sequence of times (though might have a continuous time series instead).
We now describe the case where the manifold $M$ is the circle $S^1$.

{\bf ``Projections'' of a torus to a circle.} 
Maps $\phi:\torus\to S^1$ have a nice representation.
Let $a = (a_1,\cdots,a_d)$ where $a_1,\cdots,a_d$ are integers and let 
$\theta = (\theta_1,\cdots,\theta_d)\in\torus$. 
The simplest $\phi$ has the form
$\phi(\theta) = a_1\theta_1+\cdots+a_d\theta_d \bmod1.$ Then $\phi$ is a continuous map of the torus to a circle.
For any initial point $\theta_0 \in \torus$, we have $\theta_n =  \theta_0 + n(a_1\rho_1+\cdots+a_d\rho_d )\bmod1$ and in this very simple case 
$ \theta_{n+1}-\theta_n = a\cdot\rho \bmod1:= a_1\rho_1+\cdots+a_d\rho_d \bmod1$ is constant and in this very special case we obtain a constant rotation rate for $\phi(\theta)$, namely

\begin{equation}\label{eqn:rho_phi}
\rho_\phi \bmod1=a\cdot\rho \bmod1.
\end{equation}
See Eq. \ref{eqn:a_dot_rho}. 
For $d=1$, Eq. \ref{eqn:rho_phi} says $\rho_{\phi}=a_1 \rho$ where $a_1$ is an integer. 
The integer $a_1$ depends on the choice of $\phi$, so even when $|a_1| = 1$ we can  
get $\rho$ for one choice and $-\rho$ for another choice.

\begin{equation}\label{eqn:rr}
\rho_\phi :=  a\cdot \rho\bmod1.
\end {equation}
We note that for every map $\phi$ of a torus to a circle, there are integers $a_j$ and a periodic function $g:\torus\to\R$ such that
\begin{equation}\label{eqn:g}
\phi(\theta) = g(\theta) +a\cdot \theta \bmod1.
\end{equation}
Computing a rotation rate for this map can be difficult. 
In fact, after we define the rotation rate below, it will turn out that Eq. \ref{eqn:rr} will still be true, independent of $g$, but this formula will not be very helpful in determining $\rho_\phi$ from the image of a trajectory, $\phi(\theta_n)$.

{\bf Changes of variables.} Define $\bar\theta := A\theta$ where $A$ is a unimodular transformation, that is, an invertible $d\times d$ matrix with integer coefficients.  
In these coordinates, Eq. \ref{eqn:plusrho} becomes
\begin{equation}\label{eqn:rhobar}
\bar\theta_{n+1} = \bar\theta_n + A\rho\bmod1.
\end{equation}
Note that $\rho$ is irrational if and only if $A\rho$ is so that concept is well defined. However, as we discuss in Section \ref{sec:rotvecset}, for a given irrational $\rho$
the set of $A\rho$ for all such matrices $A$ is dense in $\torus$. If for example we wanted to know the vector $\rho$ with 30-digit precision, every 30-digit vectors in $\torus$ would be valid approximations for an appropriate choice of coordinate matrix $A$.

We assume throughout this paper each continuous function such as those
denoted by $F, \phi,\gamma,$ and $\psi$, and each manifold is smooth,
by which we mean infinitely differentiable (denoted $C^\infty$). This 
assures rapid convergence of our numerical methods.

\textbf{Defining $\Delta$ and its lift $\hat\Delta$ for a 
projection $\phi:\torus\to S^1$ to a circle.}
Rotation rates are key characteristics of any quasiperiodic
trajectory.
Suppose 
there exists a continuous 
map $\phi:\torus\to S^1$ from the dynamical
system to a circle, 
but we only know 
the image
$\phi_n:=\phi(n\rho)$  
sequence of a trajectory $F(\theta_n)=
\theta_{n+1}=\theta_{n}+\rho\bmod1$ on a torus.
Define 
\begin{align}
\Delta(\theta) &= \phi(\theta+\rho)-\phi(\theta)\bmod1\nonumber\\
 &= g(\theta+\rho)+ a\cdot (\theta+ \rho) -[ g(\theta)+ a\cdot (\theta)]\bmod1\mbox{ (from Eq. \ref{eqn:g})}\nonumber\\
 &= a\cdot \rho + g(\theta+\rho) -g(\theta) \bmod1.\label{Eq5}
\end{align}
We say $\hat\Delta$ is a {\bf lift} of $\Delta:\torus\to S^1$ if 
(i)$\hat\Delta:\torus\to\mathbb{R}$, (ii)
$\hat\Delta$ is continuous; and (iii) $\hat\Delta(\theta)\bmod1=\Delta(\theta)$. 
Motivated by Eq. \ref{Eq5}, we define
 \begin{equation}\label{hatdelta}
   \hat\Delta(\theta) :=a\cdot \rho + g(\theta+\rho) -g(\theta).
\end{equation}
Then (i),(ii), and (iii) are satisfied so $\hat\Delta$ is a lift of 
$\Delta$.

Define $\hat\Delta_n= \hat\Delta(\theta_n)$. 
\begin{proposition}
Assume $\theta_n$ is quasiperiodic.
There is a well-defined rotation rate $\rho_\phi$,
\begin{equation}\label{eqn:delta}
\rho_\phi := \left(\lim_{N\to\infty}\frac{\sum_{n=0}^{N-1}\hat\Delta_n}{N}\right)\bmod1,
\end{equation}
and using the notation of Eq. 4,
\begin{equation}\label{eqn:a_dot_rho}
\rho_\phi = a\cdot\rho\bmod1.
\end{equation}
\end{proposition}

The existence of the limit is guaranteed by the Birkhoff Ergodic Theorem (See Theorem \ref{eqn:ergodictheorem}), which says that the limit in Eq. \ref{eqn:delta} is
\begin{align}
\int_\torus \hat\Delta(\theta)d\theta  &=
\int_\torus\big(
a\cdot \rho + g(\theta+\rho) -g(\theta)\big)d\theta\nonumber \\
&=a\cdot \rho + \int_\torus g(\theta+\rho)d\theta -\int_\torus g(\theta)d\theta\label{eq7}\\
&=a\cdot \rho,
\end{align}
since $\int_\torus d\theta = 1$ and the two integrals in Eq. \ref{eq7} are equal.
Hence $\rho_\phi =a\cdot \rho\bmod1.$\qed

Different choices of the lift $\hat\Delta$ can change $\rho_{\phi}$ by an integer, so $\rho_{\phi} \bmod1$ is independent of the choice of lift $\hat\Delta$.
The {\bf rotation rate} is this $\rho_\phi \bmod1$. 

The limit in Eq. \ref{eqn:delta} exists and is the same for all initial $\theta_0$. 

{\bf A caveat.} We note however, that in practice we do not know
$a\cdot\rho$ so in practice we need to determine numerically what $\hat\Delta$ is and we must numerically evaluate the limit.

Throughout this paper
we consider $\torus$ to be  $[0,1]^d\bmod1$, where each copy of $[0,1]$ 
is the fraction of 
revolution around a circle. Furthermore  $\theta\in
\torus$ can be treated as a set of $d$ real numbers in $[0,1)$.
This will enable us to write $x\in\R$ unambiguously as 
$$x = k+ (x \bmod1)$$ 
where $k$ is an integer.

There are cases where it is easy to compute the rotation rate $\rho_\phi$. 
If the angle always makes small positive increases, 
we can convert $\phi_{n+1}-\phi_n \bmod1$ into a small real positive number in $[0,1)$, 
and we can think of $\Delta_n=\phi_{n+1}-\phi_n$ as numbers in $(0,\alpha)$, where $0<\alpha<1$. 
The limit of the average of $\Delta_n$ is the rotation rate. 
The average of two or more angles in $S^1$ is not well defined. Hence we must average real numbers, not angles, and making that transition can be difficult.

{\bf Numerical determination of a lift $\hat\Delta$.} The essential problem in computing $\rho_\phi$ 
is the determination of a lift $\hat\Delta$ for $\phi$. Given a lift, we can compute $\rho_\phi$ using
Eq. \ref{eqn:delta}.
While we know the fractional part of $\hat\Delta$ is $\Delta\in [0,1),$ as we will explain later,
we must choose the integer part $k_n$ of each $\hat\Delta_n$ so that all of the points $(\theta_n, \hat\Delta_n) := (\theta_n, k_n+\Delta_n)$ 
lie on a connected curve in $S^1\times\R$ (for $d = 1)$ or a connected surface in in $\torus\times\R$ (for $d > 1)$. 
We must choose these integer parts despite the fact that we do not know which $\theta_n$ corresponds to $\Delta_n$.

\begin{figure}
\begin{center}
\includegraphics[width=.44\textwidth]{\Path 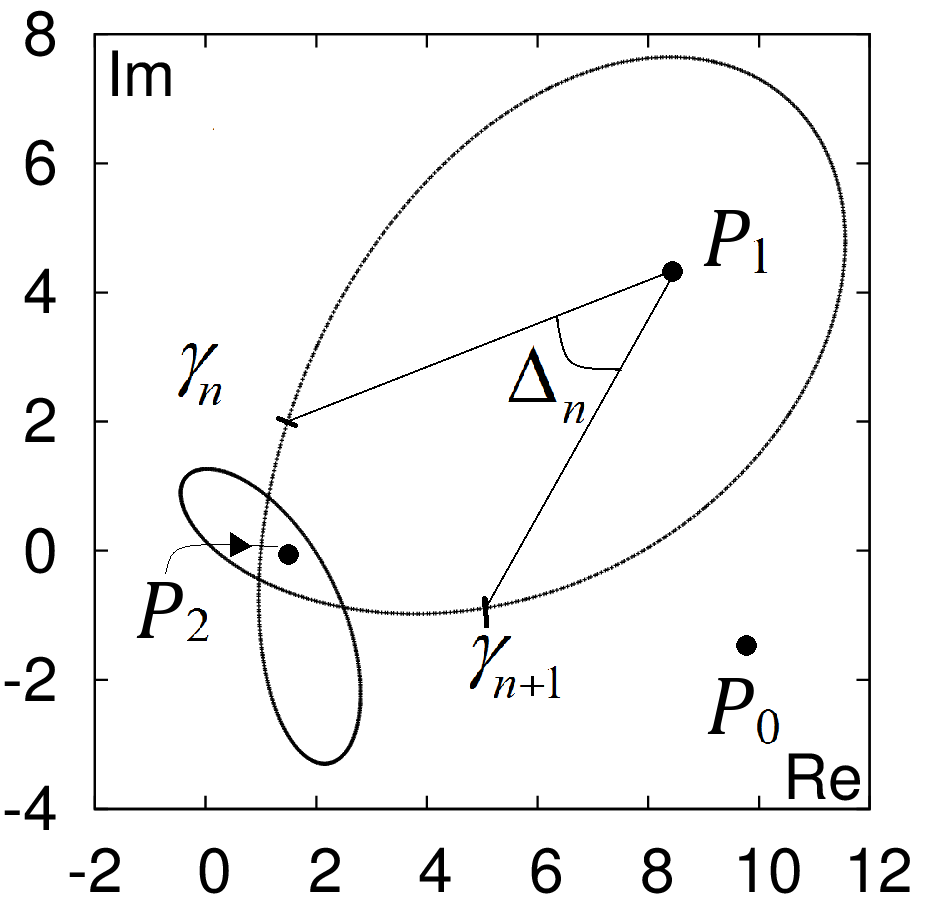}
\includegraphics[width=.42\textwidth]{\Path 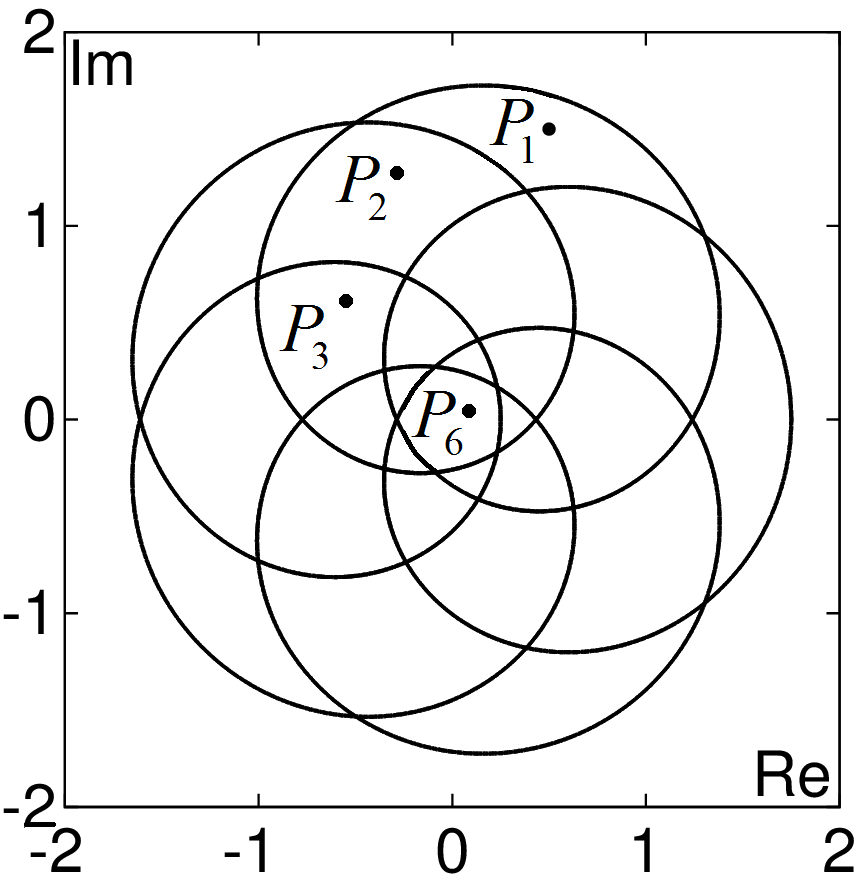}
\caption{
{\bf The fish map (left) and flower map (right).} 
The function $\gamma:S^1\to \mathbb{R}^2$ for each panel is respectively Eq. \ref{eqn:fish} and Eq. \ref{eqn:flower} and the image plotted is $\gamma(S^1)$. These are images of quasiperiodic curves with self-intersections, and we want to compute the rotation rate only from knowledge of a trajectory $\gamma_n\in \R^2$. The curves winds $j$ times around points $P_j$, so $P_1$ is a correct choice of reference point from which angles can be measured to compute a rotation rate. If instead we choose $j\ne 1$, then the measured rotation rate will be $j$ times as big as for $j=1$. In both cases, $P_1$ is the reference point. $P_1=(8.25,4.4)$ and $(0.5,1.5)$ for the fish map and flower map, respectively. The angle marked $\Delta_n\in [0,1)$ measured from point $P_1$ is the angle between trajectory points $\gamma_n$ and $\gamma_{n+1}$.}
\label{fig:fish_flower} 
\end{center}
\end{figure}

\begin{figure}
\begin{center}
\includegraphics[width=.4\textwidth]{\Path 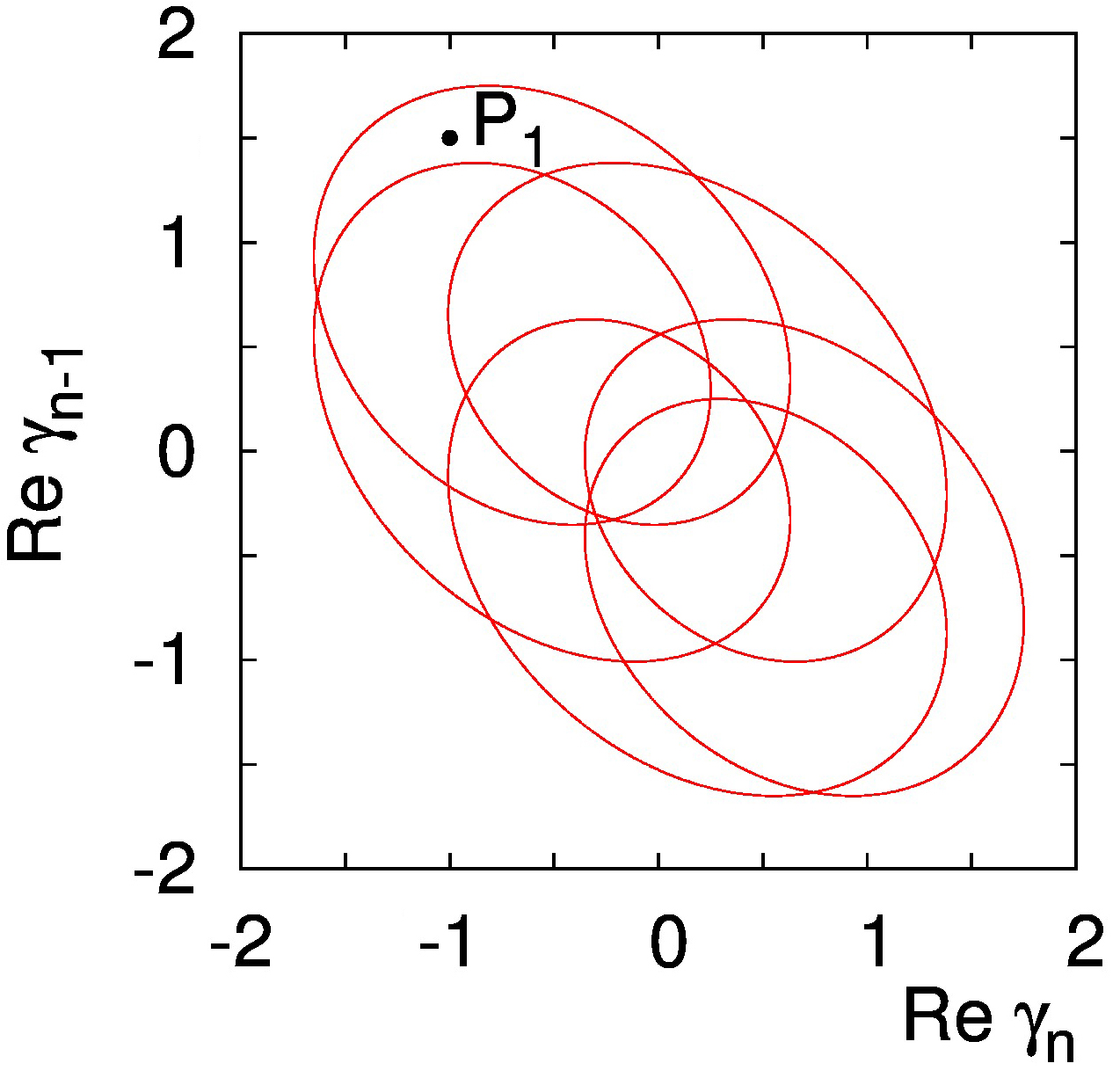}
\caption{
{\bf The flower map revisited.} 
Suppose instead of having the function $\gamma:S^1\to \mathbb{R}^2$ for the flower Eq. \ref{eqn:flower} in Fig. \ref{fig:fish_flower}, we had only one coordinate of $\gamma$, for example, the real component, $Re~\gamma.$ Knowing only one coordinate would seem to be a huge handicap to measuring a rotation rate. But it is not. In the spirit of Takens's idea of delay coordinate embeddings explained in detail later, we plot $(Re~\gamma_n,Re~\gamma_{n-1})$ and choose a point $P_1$ as before, and the map is now two dimensional. The rotation rate can be computed as before. The rotation rate $\rho_{\phi}$ here using $P_1$ is the same as for Fig. \ref{fig:fish_flower} right. 
}
\label{fig:flower_delayed} 
\end{center}
\end{figure}

\begin{figure}
\includegraphics[width=.45\textwidth]{\Path 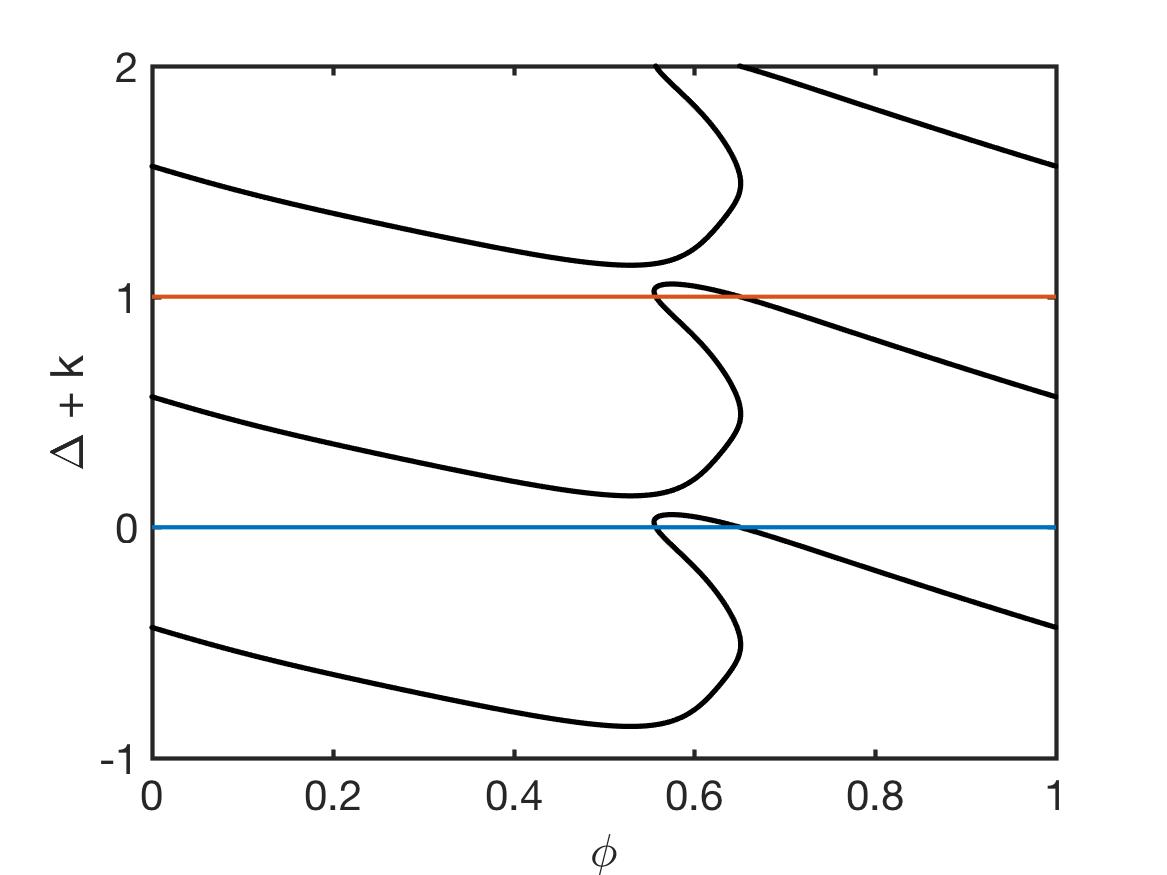}
\includegraphics[width=.45\textwidth]{\Path 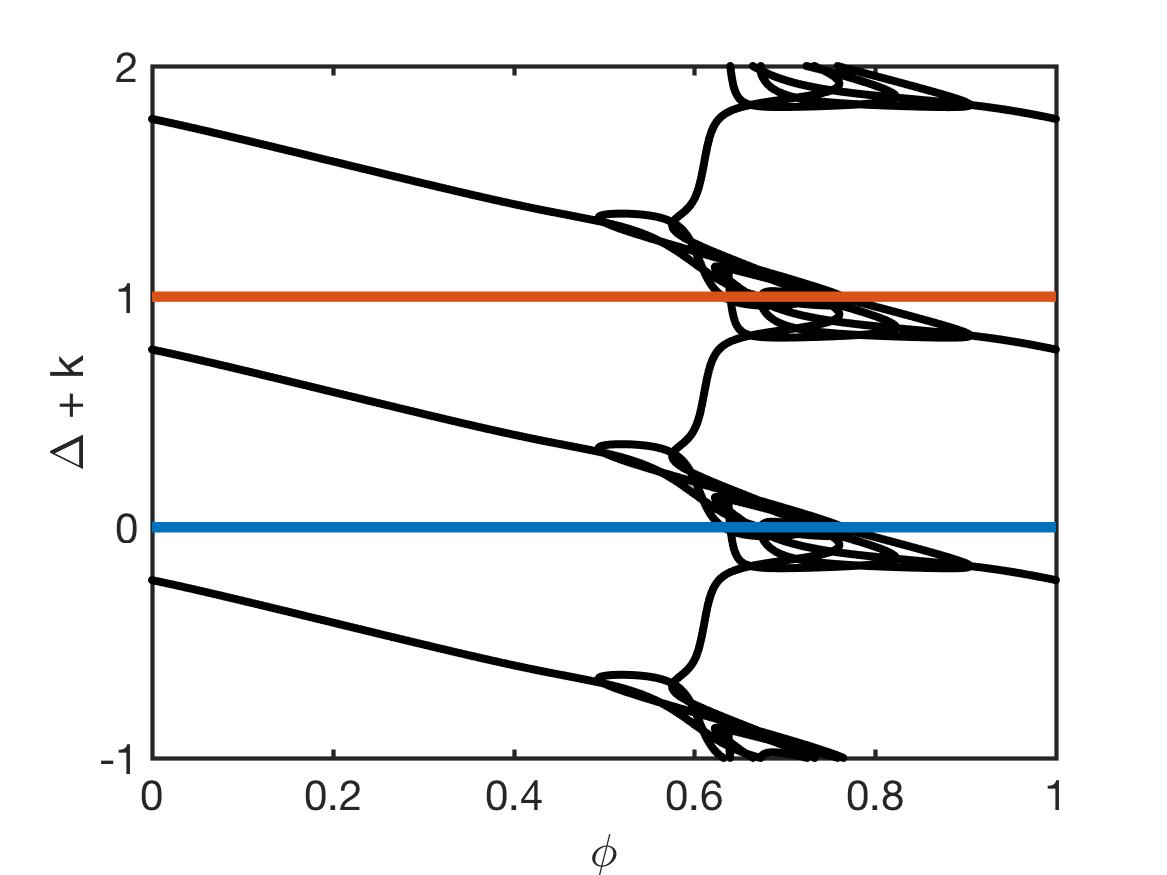}
\caption{{\bf The angle difference for the fish and the flower maps. } Here we plot  $(\phi_n, \Delta_n+ k)$ for every $n\in\mathbb{N}$ and all integers $k$, where $\Delta_n=\phi_{n+1}-\phi_n \bmod1$.  In the left panel (the fish map, the easy case) the closure of the figure resolves into disjoint sets (which are curves $\subset \mathbb{R}\times S^1$), while on the right (the flower map, the hard case) they do not. Hence if we choose a point plotted on the left panel, it lies on a unique connected curve that we can designate as $C\subset S^1\times \mathbb{R}$. We can choose any such curve to define $\hat\Delta_n$, namely we define $\hat\Delta_n = \Delta_n + k$ where $k$ is the unique integer for which $(\phi_n,\Delta_n + k)\in C$.   
A better method is needed to separate the set in the right panel into disjoint curves -- and that is our embedding method. }
\label{fig:flw_fish_1D_unlift}
\end{figure}

\begin{figure}
\includegraphics[width=.45\textwidth]{\Path 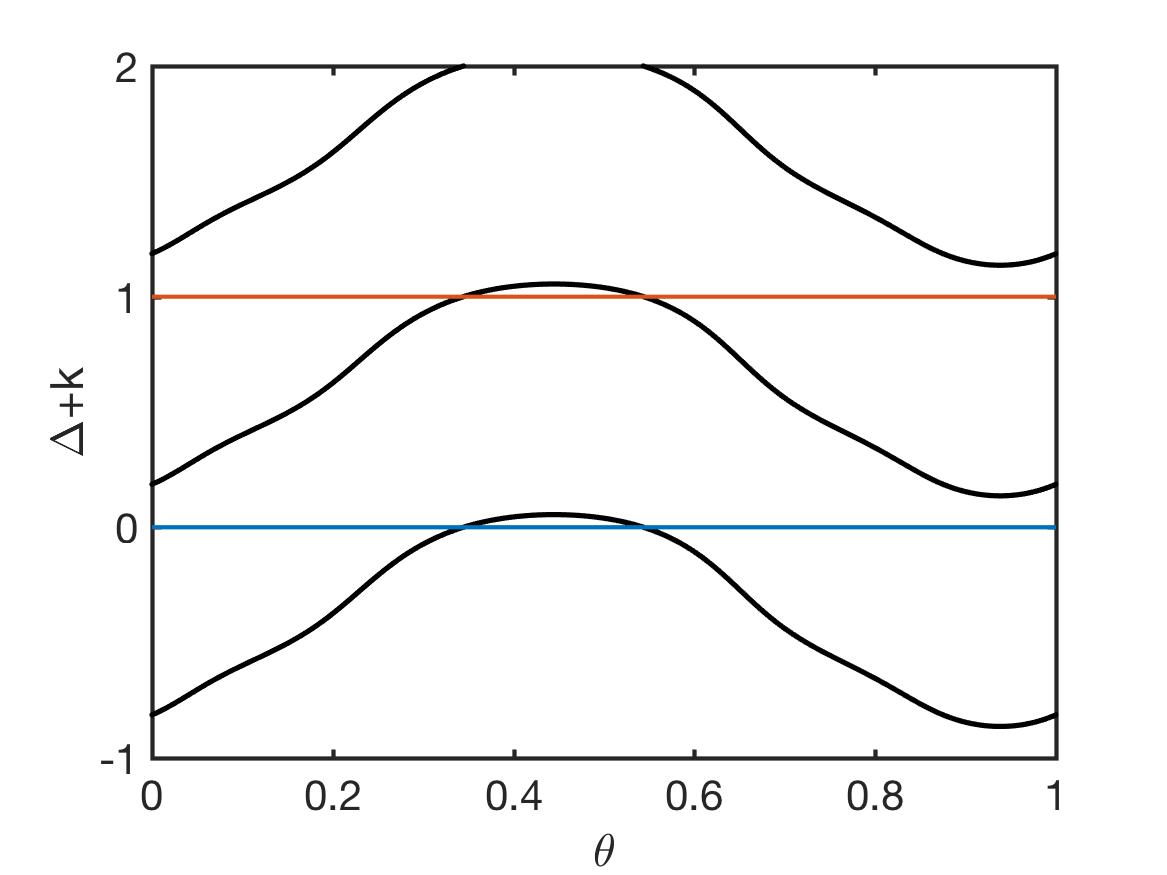}
\includegraphics[width=.45\textwidth]{\Path 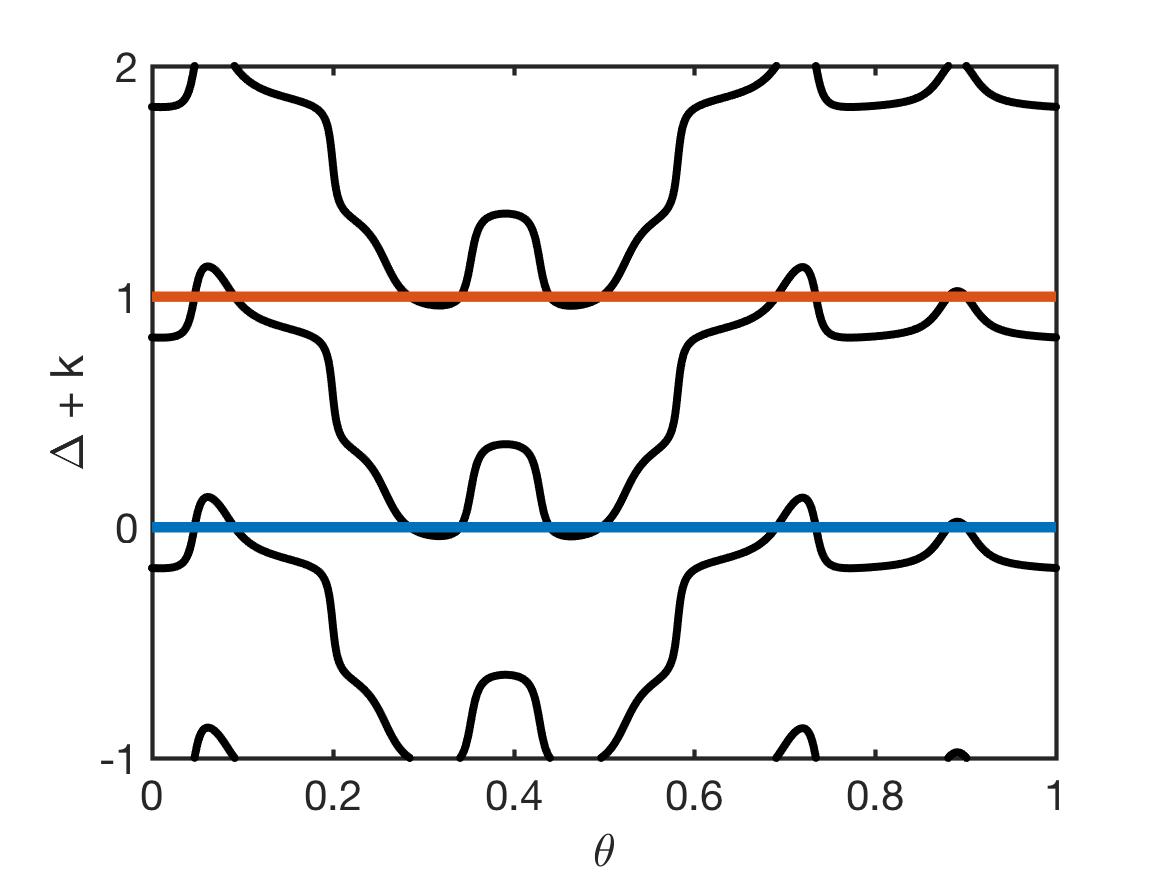}
\caption{
{\bf A lift of the angle difference for the fish and for the flower maps.}
This is similar to Fig. \ref{fig:flw_fish_1D_unlift} except that the horizontal axis is $\theta$ instead of $\phi$. 
That is, 
we take $\theta_n$ to be $n\rho$ and
$\Delta(\theta) = \phi(\theta+\rho)-\phi(\theta) \bmod1 \in [0,1)$
and we plot
$(\theta_n,\Delta_n + k)$ for all integers $k$ (where again $\Delta_n=\Delta(n\rho)$), 
These are points on the set
$G = \{(\theta, \Delta(\theta)+ k):\theta\in S^1, k\in\mathbb Z \}$.
This set $G$ consists of a countable set of disjoint compact connected sets, ``connected components'', each of which is a vertical translate by an integer of every other component. 
For each $\theta\in S^1$ and $k\in \mathbb Z$ there is exactly one point $y\in [k,k+1)$
for which $\theta,y)\in G$.
Each connected component of $G$ is an acceptable candidate for $\hat\Delta$.
Unlike the plots in Fig. \ref{fig:flw_fish_1D_unlift}, 
$G$ always splits into disjoint curves.
Unfortunately the available data, the sequence $(\phi_n)$ only lets us make plots like Fig. \ref{fig:flw_fish_1D_unlift}. 
But the Takens Embedding method allows us to plot something like $G$ and determine the lift in the next figure. 
}
\label{fig:flw_fish_1D_lift}
\end{figure}

\begin{figure}
\includegraphics[width=.80\textwidth]{\Path 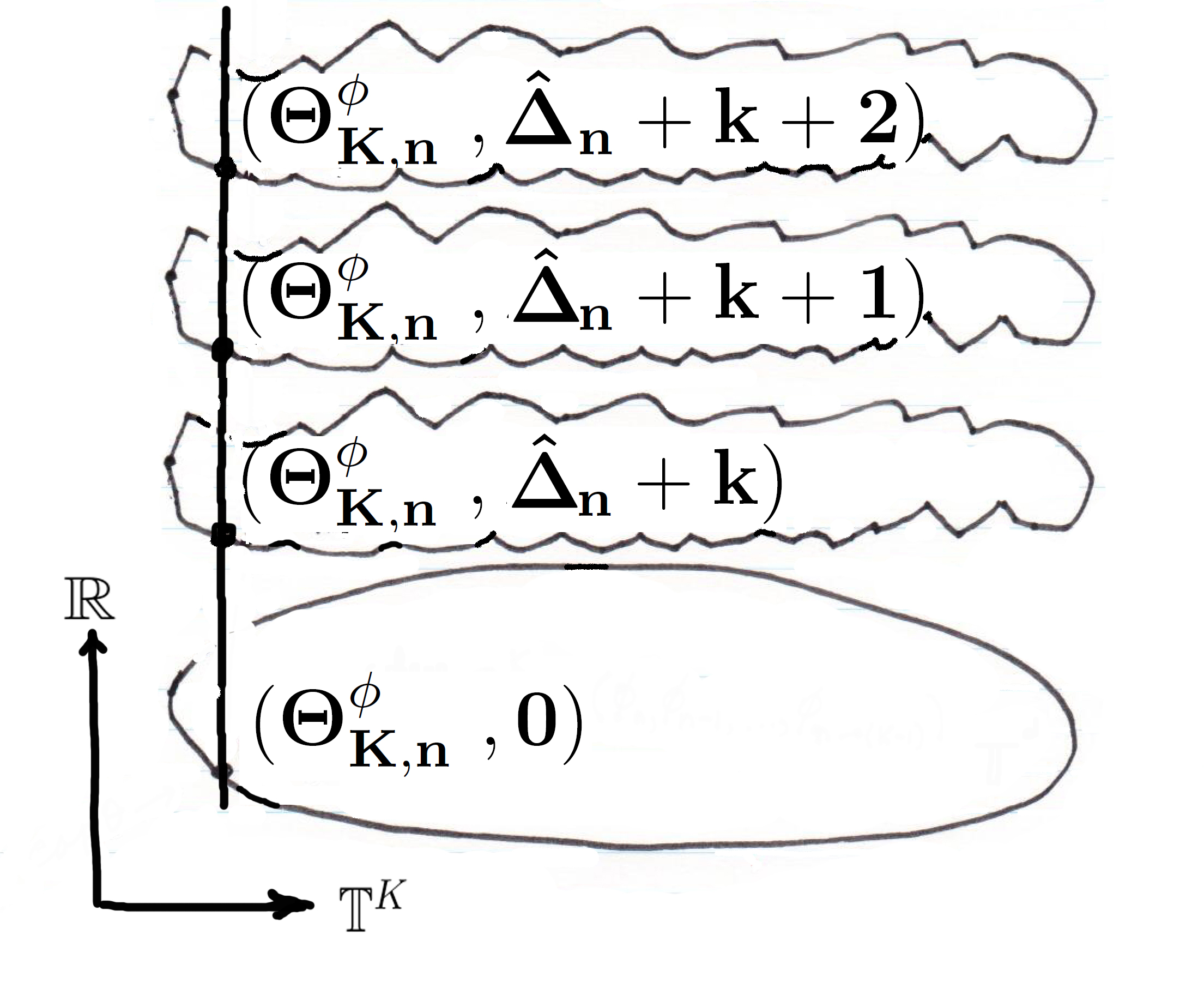}
\caption{{\bf Lifts over an embedded torus.} Let $\Theta := \Theta_K^\phi$ be as in Eq. \ref{eqn:def:embedding} and let $\theta_n = n\rho$ be a trajectory on $\torus$. Assume $K\ge 3$. By Theorem \ref{thm:takens}  for almost any map $\phi$, the set $\Theta(\torus)$ is an embedding of $\torus$ into $\TTT^{K}$; i.e., $\Theta$ is a homeomorphism of $\torus$ (the circle $S^1$ when $d=1$) onto $\Theta(\torus)$. 
In particular the map is one-to-one. The smooth (oval) curve is the set 
$(\Theta(\torus),0)$. 
As in our previous graphs, the vertical axis shows the angle difference  $\Delta ( \theta ) \in [0,1)+k$ for all integers $k$. 
Write $\mathbb{U} :=\{(\Theta(\theta),\Delta(\theta)+k):\theta\in\torus \mbox{ and } k\in\mathbb{Z}\}$. 
Unlike Fig. \ref{fig:flw_fish_1D_unlift} but like Fig. \ref{fig:flw_fish_1D_lift}, $\mathbb{U}$ always splits into bounded, connected component manifolds that are disjoint from each other.
Hence $\mathbb{U}$, which is also the closure of the set $\{(\Theta(\theta_n),\Delta_n + k):k\in \mathbb Z,n=0,\cdots,\infty\}$, separates into disjoint components each of which is a lift of $\Delta$ and each of which is homeomorphic to $\torus$. For each integer $k$ the set $\{(\Theta(\theta),\Delta(\theta)+k):\theta\in\torus\}$ is a component as shown in this figure. See Theorem \ref{thm:ourresult}.}
\label{fig:sketch-over-circle}
\end{figure}

Even in that case $d=1$ there has been no general method for computing
the lift in order to find $\rho_\phi$, though there is a literature
dealing with special cases. See for example
\cite{belova:16,luque:villanueva:14,luque:Rot}.
We have established a general method for determining the lift $\hat\Delta$, 
as summarized in
the Figs. \ref{fig:fish_flower}-\ref{fig:sketch-over-circle}. Our method is
based on the Theorem \ref{thm:takens}, a version of the  Embedding
Theorems of Whitney and Takens, described in detail in Section \ref{sec:method}.

{\bf Defining $\phi$ from a planar projection $\gamma$.} 
Assume that we are given a planar projection $\gamma:\torus\to\R^2$ and the 
images $\gamma(\theta_n)$. 
Fix a {\bf reference point} $P \in \R^2$ that is not in the image 
$\gamma(\torus)$. Let $\R^2$ be the complex plane $\mathbb{C}$,
so that we can define
$\phi(\theta)\in [0,1)\bmod1 = S^1$ by 
\begin{equation} 
e^{2 \pi i \phi(\theta)} = \cfrac{\gamma(\theta) -P}
{\| \gamma(\theta) -P \|}. 
\label{eqn:phidef} 
\end{equation}

The {\bf winding number} around $P$ is 
\begin{displaymath}
W(P) := \int_0^1  \phi' (\theta +s) \;ds, 
\end{displaymath}
where $\phi' = \frac{d\phi}{dt}$. Note that $W(P)$ is an integral over
the circle so it does not depend on $\theta$. The value of $W$ is
piecewise constant and integer-valued. In our examples, it is critical
that the projection of our quasiperiodic trajectory into $\R^2$ is such
that there exists a point $P$ in $\R^2$ with $|W(P)|=1$. 
That is because the measured rotation rate will be higher by a factor of 
$|W(P)|$. For degenerate cases, there may be no point for which $|W(P)|=1$, as shown in the next paragraph.

{\bf A non-generic map $\gamma$.} 
Consider the
map given by $\gamma(z)=z^2$ where $z \in \mathbb{C}$.  The map
$\gamma$ maps the unit circle onto the unit circle and for any value of $P \in \mathbb{C}$, 
$W(P)=0$ if the
reference point $P$ is outside that circle, and $W(P)=2$ if inside, and $W(P)$ is 
not defined if $P$ is on the unit circle. Thus there is no point $P$ such that 
$W(P) = 1$.

{\bf Two illustrative examples of complicated images of a quasiperiodic process.} Figure~\ref{fig:fish_flower} shows the
projections  maps $\gamma:S^1\to \mathbb{R}^2$, showing how the winding
number differs in different connected components of the figure. On the
left panel, every point inside the interior connected region that
contains $P_1$ can act as a reference point for measuring angles and
yields the same value of $\rho_\phi$. If the map is sufficiently simple, (i.e., the nonlinearity $g$ in Eq. \ref{eqn:g} is sufficiently small),
the rotation rate can immediately be computed as the average of these
angle differences. However, if the map $\gamma$ is more complicated,
measurement of angle is compounded by  overlap of lifts of the angle
between two iterates, since they can be represented by multiple values 
(values differing
by an integer). 

{\bf Projections to $\R$.} Sometimes we are only provided with a
scalar-valued function $\gamma: \torus \to \R$, and yet we can
still construct a two-dimensional map and use the methods described for $\R^2$ projections. 
For example,  Fig.
\ref{fig:flower_delayed} shows how we can recover a planar map from only the 
first component $Re~\gamma_n$ 
of the flower map by considering planar points $(Re~\gamma_{n-1}, Re~\gamma_{n})$. 
This map still gives  same rotation
rate as obtained by using the map in Fig. \ref{fig:fish_flower}. 

A similar example occurs with the Moon. The mean time between lunar
apogees is 27.53 days, slightly longer than the 27.3-day sidereal month.
Suppose we measure the distance $D_n$ between the centers of the Earth
and Moon once each sidereal month, $n=0,1,2,\cdots$. Then the sequence
$D_n$ has an oscillation period of 8.85 years and can be measured using
our approach by plotting $D_{n-1}$ against $D_n$, and the point
$(D_{n-1},D_n)$ oscillates around a point $P= (D_{av},D_{av}),$ where
$D_{av}$ is the average of the values $D_n$. Small changes in $P$ have
no effect on the rotation rate. 

Yet another case arises from The Moon's orbit being tilted about 5
degrees from the Earth-Sun plane. The line of intersection where the
Moon's orbit crosses the  Earth-Sun plane precesses with a period of
18.6 years. The plane of the ecliptic is a path in the distant stars
through which the planets travel. Measuring the Moon's angular distance
from this plane once each sidereal month gives scalar time series with
that period of 18.6 years. This example can be handled like the apogee
example above.

As a last example, see also our treatment of the circular planar restricted three body
problem in Section \ref{sec:3BP}  where we compute two rotation rates
of the lunar orbit, the first by plotting the rotation rate around a
central point and the second by plotting $(r,dr/dt)$, deriving the
rotation rate from a single variable $r(t)$, the distance from a central
point, where $t$ is time.

{\bf The Birkhoff Ergodic Theorem.} 
This theorem assumes there is an invariant set, which in our case is the set $\torus$. Since we are interested here only in quasiperiodic dynamics, we can assume the dynamics are given by Eq. \Xx\ where $\rho$ is irrational. 
Lebesgue measure is invariant; that is, each measurable set $E$ has the same measure as $F(E) = E+\rho$ and as $F^{-1}(E) = E-\rho$.
This map is ``ergodic'' because if $E$ is a set for which $E = F(E) = E+\rho$,
then the measure of $E$ is either $0$ or $1$. 

The measure $\mu$ enables the computation of the space-average $\int_\torus fd\mu$ for any $L^1$ function $f:\torus\to\R$ when a time series is the only information available. Since $\mu$ is Lebesgue measure, we can rewrite that integral as $\int_\torus f(\theta)d\theta$. 
We note that the Lebesgue measure of the entire torus is 1, so Lebesgue measure is a probability measure. Hence
$\int_\torus d\theta=1$.

For a map $F: \torus \to \torus$, the {\bf Birkhoff average} of
a function $f:\torus\to\R$ along the trajectory $\theta_n = F^n \theta_0$ is
\begin{equation}\label{eqn:B}
B_N(f)(\theta_0) : = \frac{1}{N}\sum_{n=0}^{N-1} f(\theta_n).
\end{equation}
\begin{theorem}[Quasiperiodic case of the Birkhoff Ergodic Theorem \cite{BrinStuck}]\label{eqn:ergodictheorem} 
Let $F: \torus \to \torus$ satisfy Eq. \ref {eqn:plusrho} where $\rho\in\torus$ is irrational. Let $\mu$ be Lebesgue measure on $\torus$.  Then 
for every\footnote{The ergodic theorem for general ergodic maps replaces ``for every'' with ``for almost every'' but for quasiperiodic maps the 
``almost'' can be omitted.} initial
$\theta_0\in\torus$, 
$\lim_{N\to\infty}B_N(f)(\theta_0)$ exists and equals $\int fd\mu$.
\end{theorem}

{\bf The Weighted Birkhoff Averaging method} ($\Q^{[p]}_N$). 
We have recently developed a method for speeding up the convergence of 
the Birkhoff sum in Theorem \ref{eqn:ergodictheorem} through introducing a $C^{\infty}$ weighting function by orders of magnitude when the process is
quasiperiodic and the function $f$ is $C^\infty$, a method we describe in ~\cite{DSSY,EPL,Das-Yorke}. 
In \cite{Das-Yorke} it is proved that the limit of using $\Q^{[p]}_N$ is the same as Birkhoff's limit.

Weighted Birkhoff ($\Q^{[p]}_N$) average of $f$ is calculated by  
\begin{equation}\label{eqn:QN}
\Q^{[p]}_N(f)(\theta_0) :=\sum_{n=0}^{N-1} \hat{w}^{[p]}_{n,N}f(\theta_n),\mbox{ where }\hat{w}^{[p]}_{n,N}=\frac{w^{[p]}(n/N)}{\sum_{j=0}^{N-1}w^{[p]}(j/N)}, 
\end{equation}
where the $C^{\infty}$ weighting function $w$ is chosen as 
\begin{equation}\label{eqn:weight}
w^{[p]}(t) :=\begin{cases}
\exp\left(\cfrac{-1}{t^p(1-t)^p}\right), & \mbox{for } t\in(0,1)\\
0, & \mbox{for } t\notin(0,1).
\end{cases}
\end{equation} 

In our calculations of the rotation rates, we use $p=1$ or $2$. 
See in particular \cite{DSSY} for details and a discussion of how the method relates to other approaches.
Note that essentially the same weight function as for $p=1$ case is discussed by Laskar~\cite{Laskar99} in the Remark~2 of the Annex, but he does not implement it.

{\bf Delay Coordinate Embeddings.} For manifolds $M_1$ and $M_2$, a map $h:M_1\to M_2$ is an {\bf embedding} (of $M_1$) if $h$ is a diffeomorphism of $M_1$ onto its image $h(M_1)$. In particular the map must be one-to-one.

Let $\psi:\torus\to M_0$ be $C^2$ where $M_0$ is a smooth manifold of dimension $D$.
In our applications below, $\psi$ is either $\phi:\torus\to S^1$ or $\gamma:\torus\to\R^2$.
While $d$ is the dimension of the domain $\torus$ of $\psi$,
$D$ is the dimension of the range. 

For a positive integer $K$, define 
$\Theta_{K}^\psi:\torus\to(M_0)^K$ as
\begin{equation}\label{eqn:def:embedding}
\Theta(\theta) :=\Theta_{K}^\psi(\theta):=\bigg (\psi(\theta), \psi(F(\theta)),\cdots, \psi(F^{K-1}(\theta))\bigg )\mbox{ for }\theta\in\torus.
\end{equation}
$K$ is referred to as {\bf the delay number} and is more precisely the number of coordinates used in defining $\Theta$. See Discussion, Section \ref{discussion}.
In the theorem below, if $K=1$, we have a Whitney-type embedding theorem, or if $D=1$, a Takens-like result. 

In order to include both of the projection maps $\phi:\torus\to S^1$ and $\gamma:\torus\to\R^2$, we introduce the more general notation $\psi:\torus\to M_0$, where the manifold $M_0$ is $D$-dimensional. 
Hence $\phi$ or $\gamma$ can be substituted for $\psi$ with $D=1$ or $2$, respectively.
\begin{theorem}\label{thm:takens}
[Special case of Theorem 2.5 in \cite{{Embedology}}]
Let $M_0$ be a smooth $D$-dimensional manifold. 
Assume $F:\torus\to\torus$ is quasiperiodic where $F$ is given in Eq. \ref{eqn:plusrho} and $\rho$ is irrational. Assume 
$$2d+1 \le KD.$$
Then for almost every $C^2$ function $\psi:\torus\to M_0$, the map $\Theta:\torus\to M_0^K$ is an embedding of $\torus$.
\end{theorem}

 While this result gives a lower bound on the delay number $K$, it is often convenient to choose $K$ much larger than required.

Define $\Gamma = \Gamma_{K}^\psi:\torus\to(M_0)^K\times \R$ as 
\begin{equation}\label{eqn:Gamma}
\Gamma(\theta) := \Gamma_{K}^\psi(\theta):= (\Theta(\theta),\hat\Delta(\theta))
\mbox{ for }\theta\in\torus.
\end{equation}
where $\hat\Delta$ is given in Eq. \ref{hatdelta}.
See Fig. \ref{fig:sketch-over-circle}.
The following corollary follows immediately from Theorem \ref{thm:takens}.

\begin{corollary}\label{cor:takens}
Assume the hypotheses of Theorem \ref{thm:takens}.
Then for almost every smooth ($C^2$) function $\psi:\torus\to M_0$, the  map $\Gamma_{K}^\psi:\torus\to M_0^K\times \R$ is an embedding of $\torus$.
\end{corollary}
Theorem \ref{thm:ourresult} explains how this result is used when we have the image of a trajectory such as $(\gamma(\theta_n))_{n=0}^{N-1}$ -- when $N$ is sufficiently large.

{\bf Comparison to previous work.} 
We have written previously about computation of rotation rate in the 
papers~\cite{DSSY,EPL,Das-Yorke}. 
A complete streamlined method for the case $d=1$ is provided in
Section~\ref{sec:method}; the Embedding continuation method is announced 
in~\cite{EPL}, but this is the first paper in which it is explained. In addition, this paper is the first time that we have applied our methods 
 to cases
where $d>1$. While we used the example (CR3BP) in \cite{DSSY}, there we used a Poincar\'e
return map  whereas here in Section \ref{sec:3BP} no return map is
used. 
We discuss the connections to our work with~\cite{belova:16, luque:villanueva:14, luque:Rot} in the subsequent sections of the paper.
Those papers do investigate the Babylonian Problem,   starting with only a set of iterates for a single finite length forward trajectory with the goal of finding a rotation number for some projection of a torus.

The investigation of quasiperiodic orbits is considered in~\cite{Haro, Cabre, Llave05, Huguet}. 
The approach in these papers assumes access to the full form of the original defining equations. 
Those papers are not investigating the Babylonian Problem.

{\bf Our paper proceeds as follows.}
We give a detailed description of our Embedding continuation method and an algorithm to implement it, in Section~\ref{sec:method}. 
Theorem \ref{thm:ourresult} gives a proof of convergence of our method.
In Section~\ref{sec:1D}, we illustrate our methods using two one-dimensional examples ($d=1$). 
We refer to these as the {\em fish map}
(introduced by 
Luque and Villanueva~\cite{luque:Rot}) and the {\em flower map}, 
based on the shapes of the graphs. 
In Section~\ref{sec:2D} we give two-dimensional ($d=2$) examples of maps for which we explore the difficulty of determining their rotation rates about a reference point.
We end in Section~\ref{sec:discussion} with a discussion.

\section{Embedding continuation method.}\label{sec:method}
We have established that there is a lift $\hat\Delta$ of $\Delta$ and that 
$\Theta$ and $\Gamma_0 :=\Gamma$ are embeddings of $\torus$ for almost every $\psi$.
{\bf We will assume in this section that $\psi$ has been chosen so that $\Theta$ and $\Gamma_0$ are embeddings.}

If we are given the image of a trajectory, either $\phi(\theta_n)$ or $\gamma(\theta_n)$, we do not yet know what the corresponding $\hat\Delta_n$ is. 
In this section, we describe how we find the lift of a map using our Embedding continuation method. 
A schematic of these ideas is depicted in Fig.~\ref{fig:sketch-over-circle}.

A major difficulty in evaluating $\rho_\phi$ is that $\hat\Delta(\theta_n)$ is not known even though $\hat\Delta(\theta) \bmod1 = \Delta(\theta)$. This is because 
$\hat\Delta(\theta)\in \mathbb{R}$ is a {\bf lift} of $\Delta(\theta)\in S^1$; i.e., they differ by an (unknown) integer $m(\theta) :=\hat\Delta(\theta) - \Delta(\theta)$. 
The key fact is that from its definition, $\hat\Delta(\theta)$ is continuous and since it is defined on a compact set it is uniformly continuous. 
We describe in Steps 1 and 2 below how to choose the integer part of $\Hat\Delta(\theta_n)$ consistently, that is, so that $\Hat\Delta(\theta_n)$ is continuous on $S^1$. 
They collectively constitute our {\bf Embedding continuation method}.

\textbf{Step 1. The embedding.}  Let $N$ be given; in practice we usually use $N\sim 10^5$ or $10^6$ if $d=1$. 
Choose the delay number $K$ so that  $2d+1 \le KD$. 
Recall that $\psi$ is either $\gamma$ or $\psi$ in our applications.
Since $\psi(\theta) \in M_0$, we have $\Theta(\theta) \in M_0^K$.
By our version of the  Takens Embedding Theorems, Theorem \ref{thm:takens}, if $2d+1 \le KD$, then
for almost every smooth function $\phi$, 
the map $\Theta$ is an embedding. 
In particular, there are no self intersections i.e., if $\Theta(\theta_1)=\Theta(\theta_2)$, then $\theta_1=\theta_2$. That implies $\Gamma$ defined by Eq. \ref{eqn:Gamma}
is also an embedding of $\torus$. We point out above that having an embedding guarantees that there are no self intersections, but there can be points far apart whose images are close to each other, and we try to avoid that by choosing $K$ large.

Denote $\mathbb{U}$ =\{$(\Theta(\theta),\Delta(\theta)+k)$ : for all $\theta\in\torus$ and all $k\in\mathbb{Z}$\}. 

{\bf The minimum distance $\epsilon$ between components of $ \mathbb{U}$.}
For each $j\in\mathbb Z$, 
define 
 $$\Gamma_j(\theta)=(\Theta(\theta),\hat\Delta(\theta)+j),$$
 and write $\Gamma_j:=\Gamma_j(\torus).$ Of course
$\Gamma_0 = \Gamma(\torus)$.
 Then $\mathbb{U}$ is the union of all $\Gamma_j$.
These sets are ``vertical'' translates of $\Gamma(\torus)$ by an integer $j$, i.e. translates in the second coordinate. 
 These are all disjoint from each other (since $\Theta(\torus)$ is assumed to be an embedding). See Fig. \ref{fig:sketch-over-circle} for an illustration.

Define 
\begin{equation}\label{eqn:epsilon}
 \epsilon := \inf\{\|p_1-p_2\|:p_1, p_2\in\mathbb U \mbox{ and are in different }\Gamma_j\},
\end{equation}
 where $\|\cdot\|$ is the Euclidean norm  on $\mathbb{R}^{2K+1}$.

 Then $\epsilon > 0$ and $\epsilon$ is the minimum distance between points on different components of $\mathbb U$. In general $\epsilon$ is hard to compute from just the time series $\psi_n:=\psi(\theta_n)$, so we have to fix a threshold $\delta>0$, assuming that $\delta<\epsilon$. 
  Then if $p_1,p_2\in \mathbb U$ and $ \|p_1-p_2\|<\delta$, it follows that  $p_1$ and $p_2$ are in the same component of $\mathbb U.$

 {\bf The choice of the delay number $K$.} It is important to note that this separation distance $\epsilon$ depends on the choice of $K$ and we observe that increasing $K$ increases $\epsilon$, so that while Theorem \ref{thm:takens} guarantees we have an embedding and therefore $\epsilon>0$, this $\epsilon$ may be small. That might make it necessary to have a very large $N$, so instead we choose $K$ much larger than the theorem requires.

\textbf{Step 2. Extending by $\delta$-continuation.} 
Write $\Theta_{n}:= \Theta(n\rho)$. The goal is to choose integers $m_n$ so that all of the points  
$(\Theta_{n},\Delta_n + m_n)$ for $n=0,\cdots,N-1$ 
are in the same component. 
This may be impossible if $N$ is not large enough.
The point $(\Theta_{0},\Delta_0)$ is in some component and we choose $m_0 = 0$ which determines a component.
Let $\mathbb A$ be the set of $n\in\{0,\cdots,N-1\}$ 
for which $m_n$ has an assigned value. This set $\mathbb A$ changes as the calculation proceeds.
Initially $m_n$ is assigned only for $n=0$ so at this point in the calculation the set $\mathbb A$ contains only $0$. 
Each time we assign a value to some $m_n$, that subscript $n$ becomes an element of $\mathbb A$.
If there is an $n_1\in \mathbb A$ and an $n_2\notin \mathbb A$ and an integer $k$ such that 
\begin{equation}\label{ineq:delta}
\|(\Theta_{n_1},\Delta_{n_1} + m_{n_1})-
(\Theta_{n_2},\Delta_{n_2} + k)\|<\delta,
\end{equation}
then the two points are in the same component and we assign $m_{n_2}=k$, which adds one element, $n_2$ to the set $\mathbb A$. Keep repeating this process (if possible) until all $\{m_n\}_{n=0}^{N-1}$ are assigned. (We will make this procedure precise in Prop. \ref{prop:steps}.)

For $N$ sufficiently large, all can be assigned values, in which case we define 
$\hat\Delta_n = \Delta_n +m_n$ for all 
$n\in\{0,\cdots,N-1\}$.
Define $$\rho_\psi^N :=
\cfrac{\sum_{n=0}^{N-1}\hat\Delta_n}{N}.$$ 

In the following theorem, we want $\delta < \epsilon$ where $\epsilon$ is in Eq. \ref{eqn:epsilon}.
\begin{theorem}\label{thm:ourresult}
For a $d$-quasiperiodic map assume $\Theta$ is an embedding. 
Given a map $\psi$, 
for $\delta$ sufficiently small, for all sufficiently large $N$ (depending on $\delta$), the above value $\rho_\psi^N$ is well defined (since all $m_n$ are defined), 
and $$\lim_{N\to\infty}\rho_\psi^N = \rho_\psi.$$
\end{theorem}


\subsection{ Continuation algorithm: long chains of little steps on $\torus$.}
To determine all $\hat\Delta(\theta_n)$ for all $n\in\{0,\cdots,n_{N-1}\}$, we begin knowing only  $\hat\Delta(\theta_0)$. Knowledge of  $\hat\Delta$ can spread like an infection, transmitted between nearby $\theta_n$. The epidemic is spread only in little steps.
The goal is to describe a continuation algorithm that identifies chains of $n_j$'s starting from $n_j=0$ and can reach every $n_j\in\{0,\cdots,n_{N-1}\}$.

To define ``little step'' we need a metric.
Let $d(\cdot,\cdot)$ be a metric on $\torus$ which 
is translation invariant, i.e. $d(x,y) = d(x+z,y+z)$ for all $x,y,z\in\torus$.
Furthermore for all $x=(x_1,\cdots,x_d)$ where all $|x_j|<0.5$, 
let $d(x,0)= \sum_j |x_j|$ (where here $d$ denotes the distance on the ``$d$''-dimensional torus).

According to Theorem \ref{thm:takens}, 
$\Theta$ is almost always an embedding of the (rigid-rotation) torus into a higher dimensional space, so we can reasonably assume the following hypothesis.

$H_1$. $\Theta$ is an embedding. (Hence $\Gamma$ is also an embedding by Cor. \ref{cor:takens}.)

In this section we will assume $\epsilon$ is given by Eq. \ref{eqn:epsilon}.
Then $(\Theta,\hat\Delta)(\torus)$ is a smooth graph over $\Theta(\torus)$. Hence if two points $\theta_1$ and $\theta_2$ in the are sufficiently close to each other, their images in $(\Theta,\hat\Delta)(\torus)$ will be less $\delta$ apart. That is 
given $\delta$ there is a $\delta_1>0$
such that 
($d(\theta_{n_1},\theta_{n_2}) < \delta_1$)
implies Ineq. \ref{ineq:delta} will be satisfied. Hence, if $m_{n_1}$ has been assigned, and $m_{n_2}$ has not, then we will now be able to assign it a value.

We say $(n_0,n_1,\cdots,n_k)$ is an $N$-$\delta_1$-{\bf chain} from $\theta_{n_0}$ to  $\theta_{n_k}$
 if $n_j\in \{0,\cdots,N-1\}$ for all $j\in\{0,\dots,k-1\}$  and $d(\theta_{n_j},\theta_{n_{j+1}}) < \delta_1$ for all $j\in\{0,\cdots,k-2\}.$

 \begin{proposition}\label{prop:steps}
 {\bf (Long Chains of Little Steps).} 
Let $F:\torus\to\torus$ be the rigid rotation with rotation vector $\rho$ with a dense trajectory. 
For $\delta_1 > 0$, there is $N>0$ such that 
for every $n\in \{0,\cdots,N-1\}$
there is a $N$-$\delta_1$-${\bf chain}$ from $\theta_{0}$ to  $\theta_{n}$. 
\end{proposition}

The following corollary interprets this proposition in terms of lifts and its proof is immediate.

 \begin{corollary}Assume $H_1$.
Assume $\delta_1>0$ is such that
$d(\theta_{n_1},\theta_{n_2}) < \delta_1$ implies Ineq. \ref{ineq:delta}.
Then, since $m_0=0$, all $m_n$ can be determined.
Write $\hat\Delta_j = \Delta_j +m_j$. Then all  $\hat\Delta_j$ are in the same lift of $\Delta$. In other words, $(\Theta^\phi_{K,j},\hat\Delta_j)$ are all in the same component of $\mathbb{U}$ 
where $\mathbb{U}$ is defined in Fig. \ref{fig:sketch-over-circle}.
\end{corollary}

To sketch a proof of the Proposition, we need the following fact.
It is an elementary fact whose proof we leave to the reader.

Given $\delta_1>0$, there exists an $N$ with the following property.

$H_2$. There exist  integers $0<\sigma_1< \sigma_2<\cdots<\sigma_P$ for some integer $P>1$ that $(i)$ the  $\sigma_j$ are relatively prime 
(i.e., the greatest common factor of all $\sigma_j$ is $1$)
and $(ii)$ $\theta_{\sigma_j}$ are within $\delta_1$ of $\theta_0$. Furthermore, $\sigma_1+\sigma_P<N$.

It is always possible to choose $N$ sufficiently large that $P=2$ in $H_2$; however, we might not want to choose such a large $N$, and we might be satisfied with having $P>2$. 

\begin{figure}
\begin{center}
\includegraphics[width=.5\textwidth]{\Path 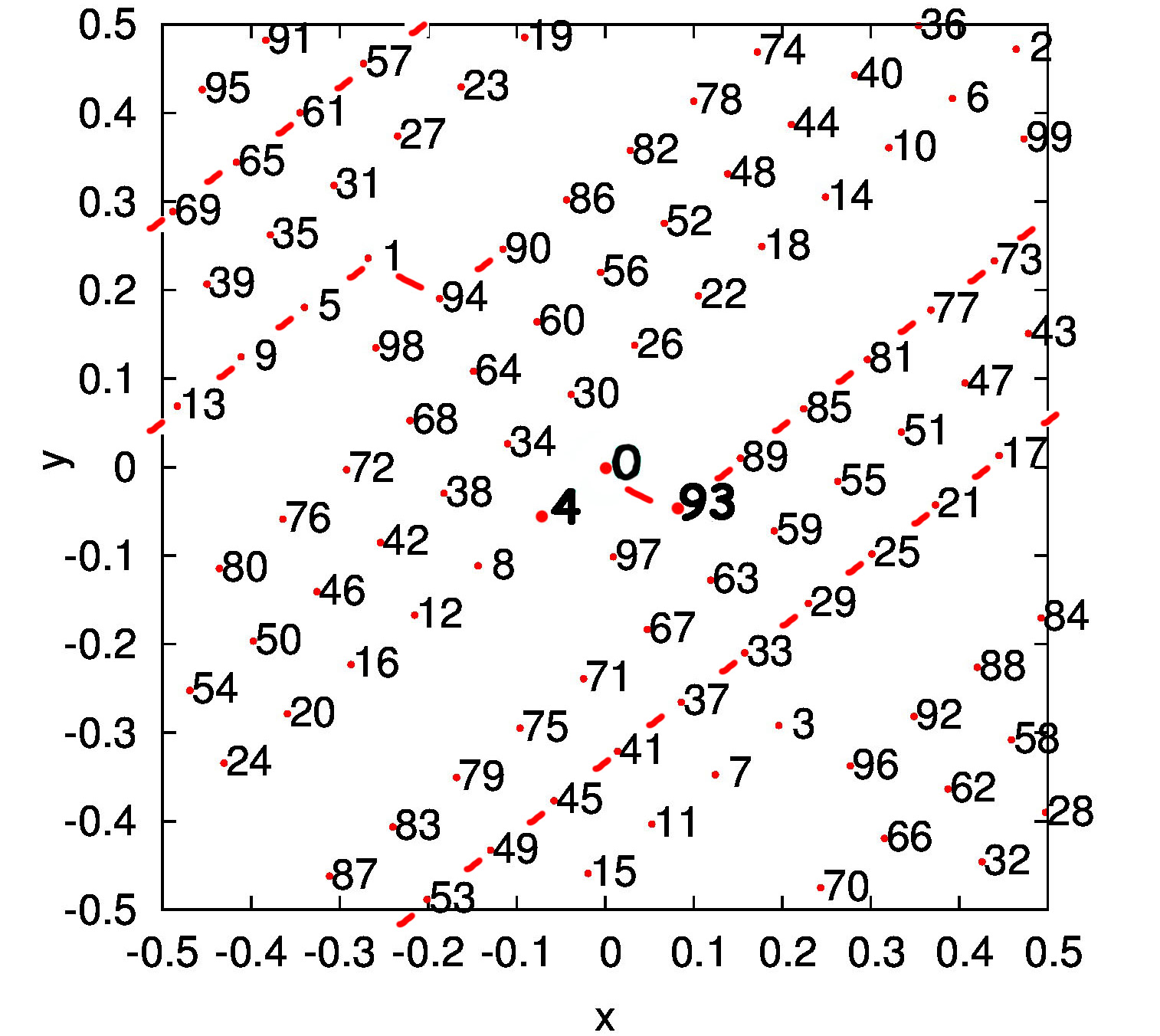}
\includegraphics[width=.47\textwidth]{\Path 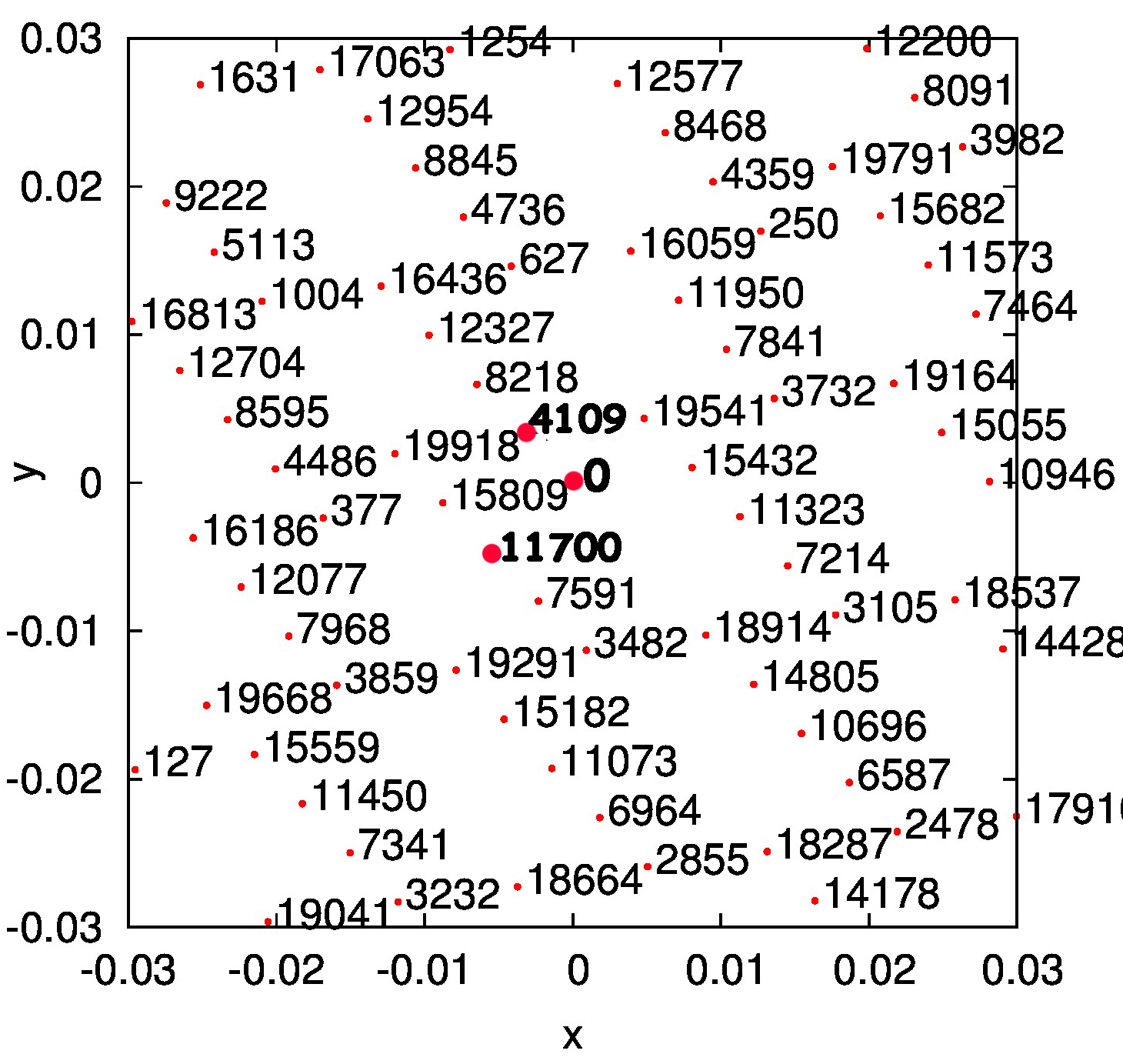}
\caption{
{\bf Rigid rotation on the torus.} 
$x_n=n\sqrt{3}~(\bmod 1),~y_n=n\sqrt{5}~(\bmod 1)$ for $n=0, \cdots, N-1$ are plotted with the origin indicated by $0$ at the center on the panel. Each point $\theta_n=(x_n,y_n)$ is labeled with its subscript $n$. 
Here $N=100$ (left) and $=20,000$ (right). Only the neighborhood of the 
origin is shown for the right panel.
In the left panel, $\theta_4$ and $\theta_{93}$ (i)  are near the origin and (ii) their subscripts are relatively prime and (iii) the total of the subscripts is less than $N$. On the right points with subscripts $4109$ and $11,700$ play the corresponding role.
In each case it follows that there is a chain of points starting from $0$ and ending at any desired $\theta_m$ where $0 < m < N$. This chain is a series of steps, each achieved by either adding one of the two subscripts or subtracting the other. See Prop. \ref{prop:steps} and the algorithm sketched in its proof. In the left panel such a chain -- adding $93$ or subtracting $4$ at each step -- is shown that ends at $\theta_{90}$.}
\label{fig:sqrt3-sqrt5} 
\end{center}
\end{figure}

{\bf An example of a pair $\theta_{\sigma_1}$ and $\theta_{\sigma_2}$ with relatively prime subscripts in dimension $d=1$.} The algorithm for creating chains does not depend on the dimension $d$.  Here we let $d=1$ and $\rho = \pi-3$ and $N= 200$ and $\delta_1 = 0.01$
(where $d(0,x)=|x|$ for $x$ close to $0$).  Then we can choose $\sigma_1=7$ and $\sigma_2= 113$ since $\pi7-22\approx 0.008$ and $355-\pi{113}\approx 0.00003$ 
so $\theta_7$ and $\theta_{113}$ are within $\delta_1$ of $0$ and the subscripts $7$ and $113$ are relatively prime. 
We can reach every subscript in $\{0,\cdots,N-1\}$ by starting from $0$ taking little steps, either increasing the subscript by $113$ or decreasing it by $7$, all the while staying between $0$ and $N$, taking steps of size less than $\delta_1$.

{\bf An example of a pair $\theta_{\sigma_1}$ and $\theta_{\sigma_2}$  with relatively prime subscripts in dimension $d=2$.}  See Fig. \ref{fig:sqrt3-sqrt5}. On the left where $N=100$, a chain is shown from $0$ taking only steps of either $+93$ or $-4$. Both are within $\delta_1=0.13$ of $0$. It would work equally well to take only steps of $-93$ or $+4$. 
When $N=20,000$ on the right, there are two relatively prime subscripts $4109$ and $11700$ whose $\theta$ values are within $\delta_1=0.011$ of $0$. 

{\bf Proof of Proposition \ref {prop:steps}.} We now describe why each $\theta_n$ can be reached by a chain starting from $\theta_0$.

We assume that for the given $\delta_1$, the $N$ and $\sigma_j$ have been chosen so that $(i)$ and $(ii)$ in $H_2$ are satisfied.

Let $B:=B(\delta_1)$ denote the $\delta_1$ neighborhood of $\theta_0=0$.
First we assume the number $P$ of $\sigma_j$ satisfies $P=2$, so 
$\theta_{\sigma_1}$ and $\theta_{\sigma_2}$ are in $B$ and their
subscripts are in $\{0,\cdots,N-1\}$ and are relatively prime.

For non-negative integers $a_1,a_2$, write 
$$[[a_1,a_2]]:=\theta_{a_2\sigma_2-a_1\sigma_1}\in \torus.$$
Suppose $a_1$ and $a_2$ are such that 
\begin{equation}\label{subscript}
0 \le a_2\sigma_2-a_1\sigma_1 < N.
\end{equation}
Since  $0 \le \sigma_1+\sigma_2<N$, we can either increase $a_1$ or $a_2$ by $1$ (thereby 
decreasing
  $a_2\sigma_2-a_1\sigma_1$ by $\sigma_1$ or increasing it by $\sigma_2$, respectively) and still have Condition \ref{subscript} satisfied.

Notice that the distance from $[[a_1,a_2]]$ to $[[a_1\pm 1,a_2]]$ or 
$[[a_1,a_2\pm1]]$ is less than $\delta_1$.
That is, changing either $a_1$ or $a_2$ by $1$ moves $[[a_1,a_2]]$ by less than $\delta_1$. 

The key step of the proof is the following.

{\bf Algorithm.} We choose a chain, a finite sequence $(\theta_j)$ of such points, each of the form  $[[a_1,a_2]]$ as follows.
Our algorithm begins at $\theta_0$ with $a_1 = a_2 = 0.$

$A_1. $ Increase $a_2$ by $1$ provided the subscript remains non negative; otherwise increase  $a_1$ by $1$. Repeat the process. Eventually the subscript returns to $0$ (with $a_1 = \sigma_2$ and $a_2 = \sigma_1$. We have thereby created a chain of points on the torus, but we most likely have not encountered all the $\theta_j$.

Next,

$A_2.$  for each point in that chain, increase $a_2$ by $1$, and repeat as long as Condition \ref{subscript} satisfied. This process yields $\theta_n$ for every $n\in \{0,\cdots,N-1\}$.

If $P=2$, we are done. 
When $P>2$, the greatest common factor, denoted $\Psi_2$, of $\sigma_1$ and $\sigma_2$ is  greater than $1$. Then the above procedure reaches all points with subscripts divisible by $\Psi_2$ and no others. The next step is essentially the same as $A_2$ except that steps are taken by adding $\sigma_3$ to 
the subscript; that is, 

$A_3$. Repeatedly add $\sigma_3$ to the subscript, as long as it remains less than $N$.

Taking all of those points and taking a small step for each by adding or subtracting $\sigma_3$ repeatedly will reach all points whose subscript is divisible by $\Psi_3 :=$ the greatest common divisor of 
 $\sigma_1$, $\sigma_2$, and $\sigma_3$. 
 
$A_j$. For each point that has been found so far,  repeatedly add $\sigma_j$ to the subscript as long as it remains less than $N$.

Eventually all $\theta_n$ for $0 <n<N$ will be reached.  \qed

\subsection{A dense set of equivalent representations for each rotation vector}
\label{sec:rotvecset} 
While the definition of quasiperiodicity requires that the map has some coordinate system that turns the map into 
Eq.~\ref{eqn:plusrho}, that requirement by itself does not determine the coordinates of $\rho$. Fixing a coordinate system allows us to write
$\rho=(\rho_1,\cdots,\rho_d)$
We have defined $ \rho $ in Eq. \ref{eqn:plusrho} in terms of a given coordinate system. 
Let $\bar\theta = A\theta$ where $\bar\theta\in\torus$ and $A$ is a unimodular transformation, that is an integer-entried matrix with determinant $|det A| = 1$, then in this new coordinate system Eq. \ref{eqn:plusrho} becomes
\begin{equation}\label{eqn:plusrho2}
\bar\theta \mapsto \bar\theta + A\rho\bmod1
\end{equation}
which is essentially Eq. \ref{eqn:rhobar}. 
Hence $A\rho$ is also a rotation vector for the same torus map. 
Below we show we have a dense set of rotation vector representations.
 
Let $\mathcal{S}$ denote the set of integer-entried $d\times d$ matrices with determinant $\pm 1$. 
Observe that for any $B\in \mathcal{S}$, $B^{-1}\in \mathcal{S}.$ A matrix in 
$\mathcal{S}$ can be viewed as a change of variables on the torus, since it preserves volume. 
Therefore we call a vector $\tilde\rho \in \R^d$ a {\bf rotation representation} of $ \rho \in \R^d $ if $\tilde\rho = A\rho $ 
for some $A\in \mathcal{S}$.
We ask: {\it When the vector $\rho$ is irrational, what are all the possible rotation
vectors (i.e., rotation representations), assuming $A\in \mathcal{S}$?} 

\begin{proposition} Assume dimension $d\ge2$. For an irrational rotation vector $\rho$, the set of its rotation vector representations is $S\rho$ (i.e., $\{A\rho :A\in S\}$), 
and $S\rho \bmod1$ is dense in $\torus.$ 
\end{proposition}

{\bf Proof.} To simplify notation we prove only the case of $d=2$. The proof for $d>2$ is analogous. See~\cite{katok:hasselblatt}. 
Write $\rho = (\rho_1,\rho_2)$. 
Note that the matrices 
$B_m :=\left(\begin{array}{c}1~~m \\0~~1 \end{array}\right)$ and $C_k :=\left(\begin{array}{c}1~~0 \\ k~~1 \end{array}\right)$ 
are in $\mathcal{S}$
for all integers $m$ and $k$, as is $A := B_m C_k$.
Then the vectors 
$(\rho_1, y_k)
=C_k (\rho_1, \rho_2)\mod~1
$ are {\em vertical} translates (translates in the direction $(0,1)$) of 
$(\rho_1, \rho_2)
\mod~1$,
where $\{y_k\} $ is a dense set in $S^1$. 
When we similarly apply $B_m$ for all $m$ to each
$(\rho_1, y_k)$
we obtain a dense set of {\em horizontal} translates of 
$(\rho_1, y_k) $
and thereby obtain a dense set in $\Torus$. 
Every coordinate of every point in that dense set is of the form 
$k_1 \rho_1 + k_2 \rho_2\mod~1$ where $k_1$ and $k_2$ are integers. \qed

\section{Examples of one-dimensional quasiperiodicity.}\label{sec:1D}

In this section, we give a detailed explanation of how we compute rotation rates for 
quasiperiodic maps on one-dimensional tori. For the
first example computation of the rotation number is easy and straight forward while in the second it is sufficiently hard that we need our method. The pair of examples makes it clear when our method should be used.

One advantage of the examples below is that we know $\rho$ and therefore we can compare it with the computed rotation rates.

{\bf Example 1. The fish map.} 
Luque and Villanueva~\cite{luque:Rot} addressed the case of a
quasiperiodic planar curve $\gamma:S^1\to \mathbb{C}$ and introduced
what we call the {\em fish map}, depicted in the left panel of Fig. \ref{fig:fish_flower}. Let
\begin{equation}\label{eqn:fish}
 \gamma(\theta):=\hat{\gamma}_{-1}z^{-1}+\hat{\gamma_{0}}+\hat{\gamma_{1}}z+\hat{\gamma_2}z^2, 
\end{equation}
where $z = z(\theta):=e^{i 2\pi \theta}$ and 
$\hat\gamma_{-1}:=1.4-2i,$ $\hat\gamma_{0}:=4.1+1.34i,$ 
$\hat\gamma_{1}:=-2+2.412i,$ 
$\hat{\gamma_2}:=-2.5-1.752i$. 
(See Fig. 5 and Eq. 31 in~\cite{luque:Rot}). They chose the rotation
rate $\rho=(\sqrt{5}-1)/2\approx 0.618$ for the trajectory $\gamma_n =
\gamma(n\rho)$ for $n = 0,1,\cdots$ so we also use that $\rho$. The
method in \cite{luque:Rot} requires a step of {\em unfolding} $\gamma$,
which our method bypasses. We measure angles with respect to $P_1 = 8.25
+ 4.4i$, where the winding number $|W(P_1)|=1$.

{\bf Example 2. The flower map.} We have created an example, the {\em flower map} in 
Fig. \ref{fig:fish_flower}, right, to be more challenging than the fish.  
Let
\begin{equation}\label{eqn:flower}
\gamma_6(\theta): = (3/4) z + z^6 \mbox{ where } z = z(\theta) := e^{i 2 \pi \theta}.
\end{equation}
We use the same $\rho=(\sqrt{5}-1)/2$ as above. The choice of a
reference point $P_1$ for which $|W(P_1)|= 1$ is shown in the right panel of 
Fig.~\ref{fig:fish_flower}. For our computations, we use $P=P_1:=
0.5 + 1.5i$. Points $P_j$ with $|W(P_j)| = j $ for $j=1,2,3,6$ 
are also shown. For $\rho_{\gamma_6}$, the rotation rate of $\gamma_6$, to yield $\rho$ or $1-\rho \bmod1$
is essential to choose a point $P$ where $|W(P)|=1$. In this example the
values of $\Delta_n$ are dense in $S^1$, and  $\max_\theta
\hat\Delta(\theta) - \min_\theta \hat\Delta(\theta) \approx 1.2$.

For both examples, Fig.~\ref{fig:fish_flower} shows two successive iterates $\gamma_n$
and $\gamma_{n+1}$, and the angle $\Delta_n$ between these two iterates,
computed with respect to a  reference point $P_1$. It was computed by finding
$\phi_n$, the angle of $\gamma_n$ with respect to $P_1$ as in
Eq. \ref{eqn:phidef}. Using this, $\Delta_n =
\phi_{n+1}-\phi_n\in[0,1)\approx S^1$. On the left, in the fish map case, if we choose $\hat\Delta_0 :=
\Delta_0$ (or alternatively  $:= \Delta_0+m$ for some $m$), then we have
selected the component in which all $\hat\Delta_n$ must lie. This is
what is referred to below as {\em the easy case}. Choose some $k$,
write $J_k := [a,b]$. Choose $m_n$ is the integer for which $\Delta_n +
m_n \in J_k$. It is not as easy to do this for the flower map on the right. 
Fig.~\ref{fig:flw_fish_1D_unlift} right shows that the possible lifts
when plotted against $\phi$ form a tangled mess which does not resolve
into bounded components, while when  plotted against $\theta$ we obtain
components that are diffeomorphic to $S^1$.

Figs. ~\ref{fig:flw_fish_1D_unlift} and ~\ref{fig:flw_fish_1D_lift}
show the possible lift values $\hat{\Delta}_n$ of the angle difference
$\Delta_n$ plotted with respect to angle $\theta$ in Fig.
~\ref{fig:flw_fish_1D_unlift} and  $\phi$ in
Fig.~\ref{fig:flw_fish_1D_lift}. For the fish map on the left, we see
that we can set $[a,b]\approx [0.18 +k,1.05+k]$ for any integer $k$.
Furthermore, we investigated the rotation rate of the signal viewed from
$P_1 = 7+4i$. Using the Weighted Birkhoff Average, 
we
observe that the deviations of the approximate rotation from $\rho$
falls below $10^{-30}$ when the iteration number exceeds $N =20,000$,
and since we know the actual rotation rate, we can report that the error
in the rotation rate is then below $10^{-30}$. Once we have found a 
proper lift for the flower map, we can do the same procedure. The next
section explains how we go about finding a lift in this more complicated case.


\section{Higher-dimensional quasiperiodic examples} \label{sec:2D}
We develop a higher-dimensional method to compute the rotation
vector $\rho$ purely from knowledge of the sequence $\theta_{n+1} :=
F(\theta_n).$ 
The question of how to compute
the rotation vector is actually two questions. Question 1: If we compute a
rotation vector, what are the possible values? Question 2: How do we compute any
of the possible values for the rotation vector in difficult cases? 
Figs~\ref{fig:flw_fish_coordinates},
 \ref{fig:flw_fish_unlift}, and
 \ref{fig:flw_fish_lift},
demonstrate that
like in one dimension, in $d$ dimensions we are able to use $d$
independent planar projections combined with a higher-dimensional
version of our Embedding continuation method in order to find a lift,
each projection leading to one component of a $d$-dimensional rotation
vector. In fact, these rotation vectors are not unique. In this section,
we give a detailed discussion of our higher-dimensional method,
describing the possible values we can achieve in calculating a rotation
vector. We then illustrate our method for three examples: the fish
torus, the flower torus, and the restricted three-body problem.

\begin{figure} \begin{center}
\includegraphics[width=.45\textwidth]{\Path 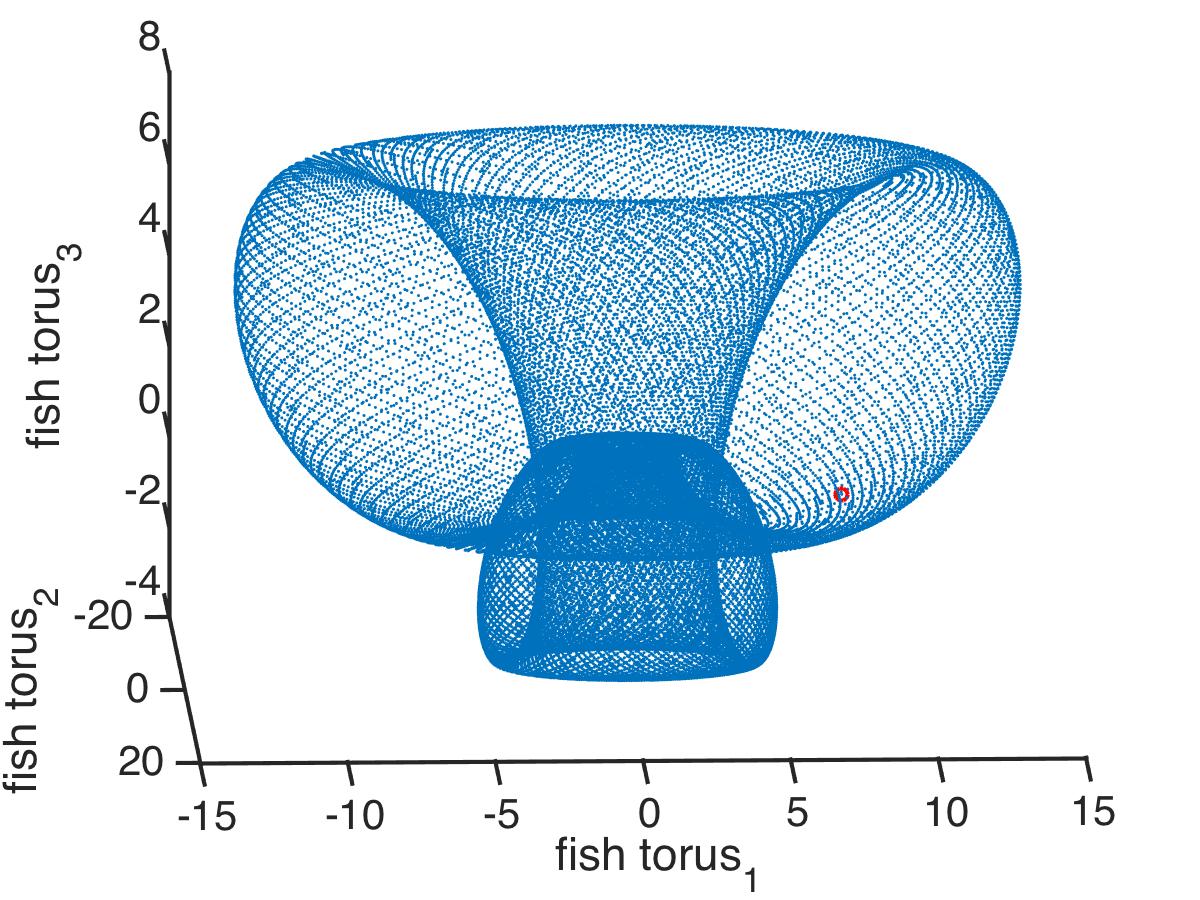}
\includegraphics[width=.45\textwidth]{\Path 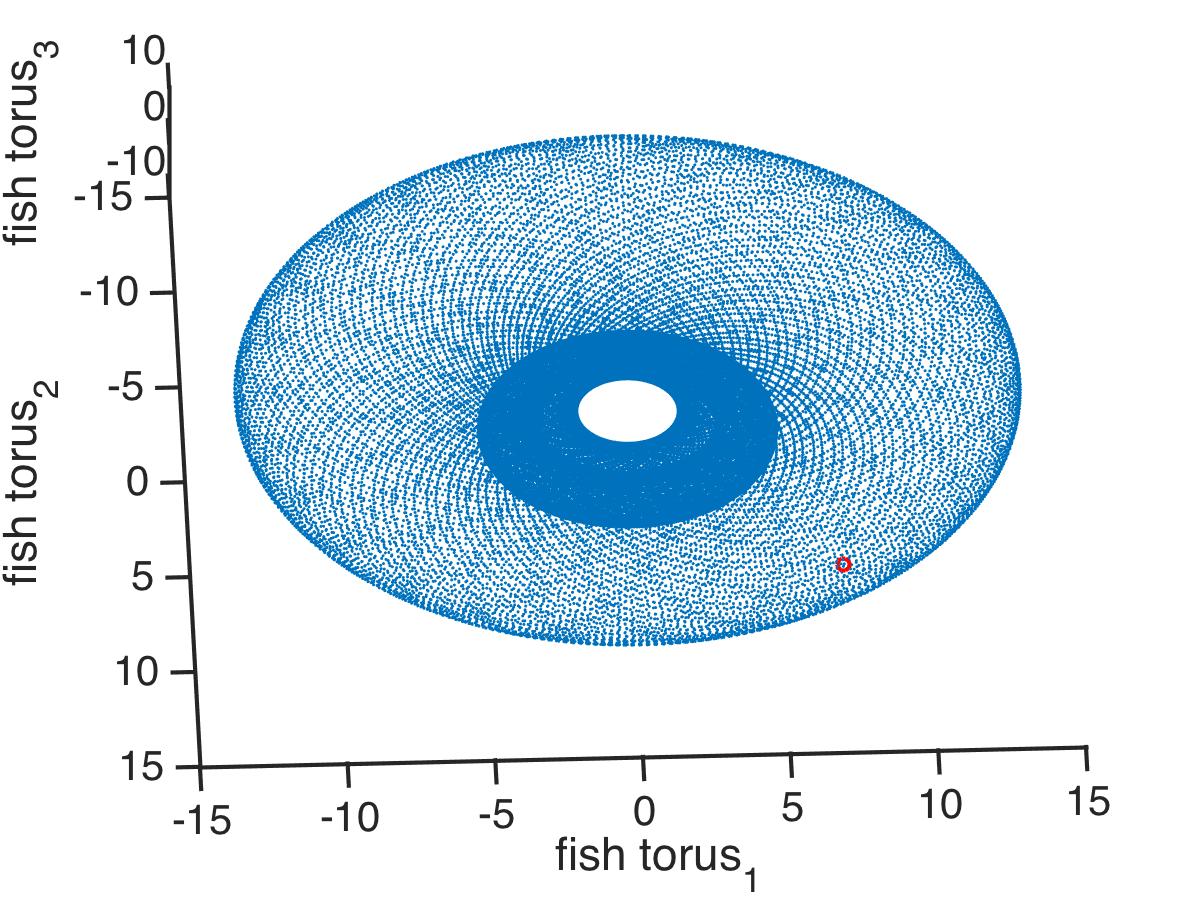}\\
\includegraphics[width=.45\textwidth]{\Path 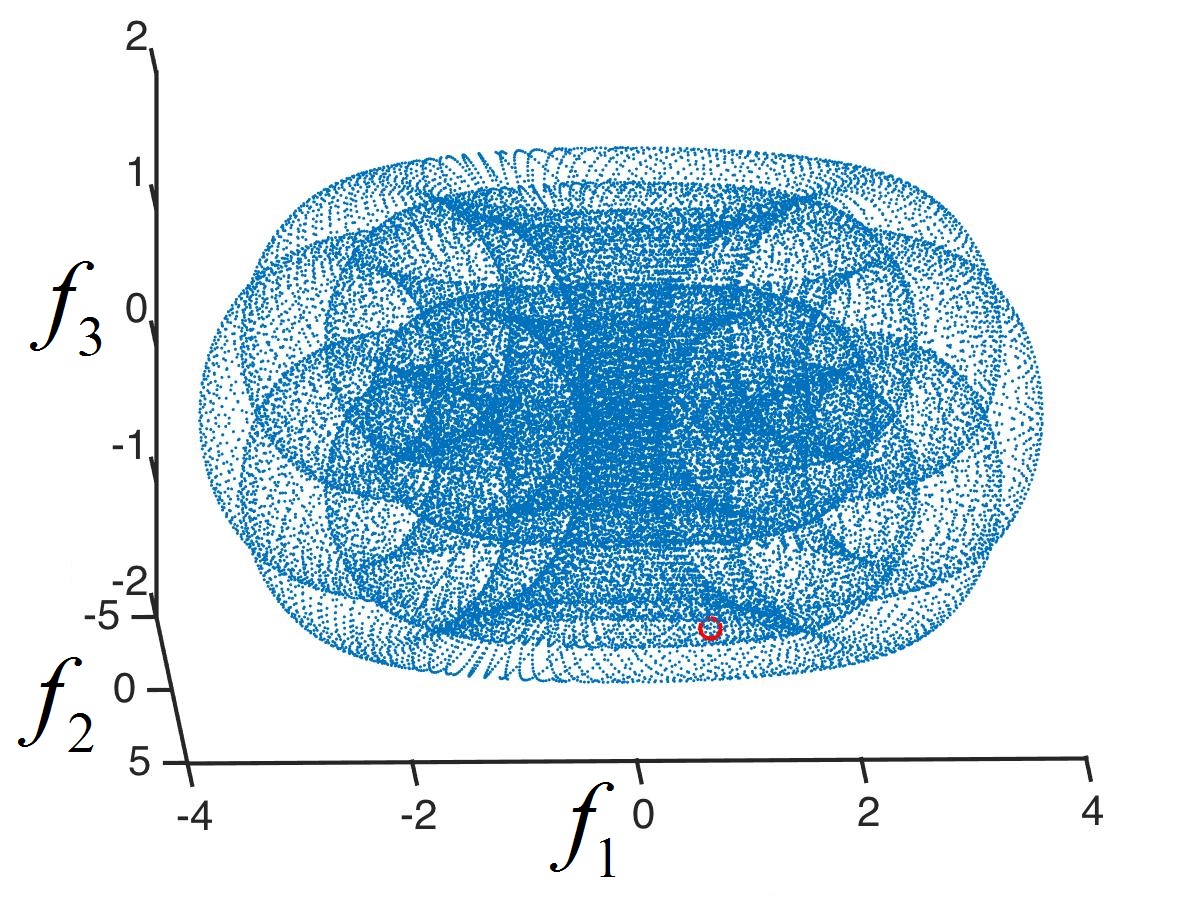}
\includegraphics[width=.45\textwidth]{\Path 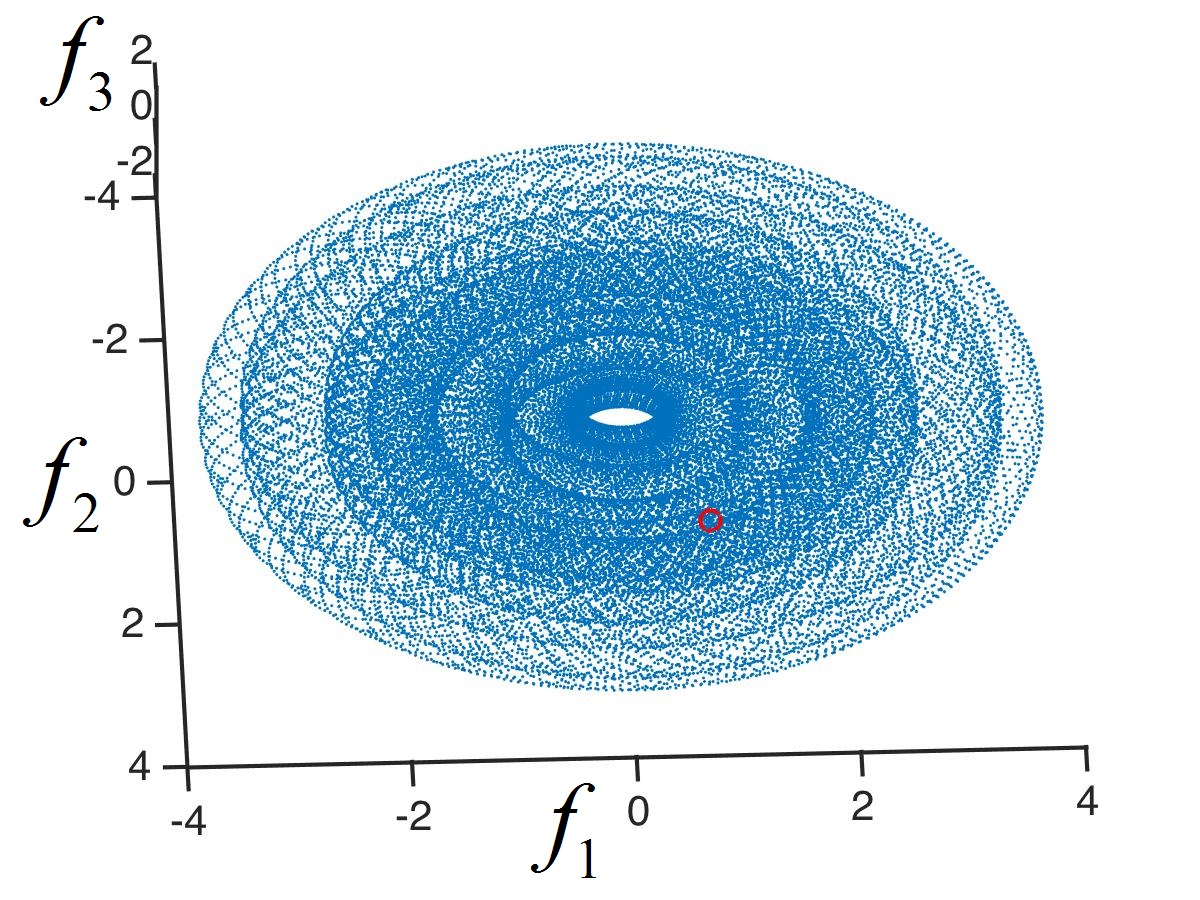}
\caption{
{\bf The fish and flower torus.} The top figures show two views of the fish torus, and 
the bottom two views of the flower torus. These figures can be thought of as projections of tori onto
the plane represented by the page. The three coordinate axes are
presented here to clarify which two-dimensional projection is being
used. The projections of the tori on the left are simply connected so
there is no way to choose a point $P$ that would yield a non-zero
rotation rate. The projections on the right 
yield images of the tori 
 that are annuli with a hole in which $P$ can be chosen to yield non-zero
results. Each is a plot of $N=50090$ iterates. The red circle is the
initial point.
}\label{fig:flw_fish_2D} \end{center} \end{figure}

\begin{figure}
\includegraphics[width=.45\textwidth]{\Path 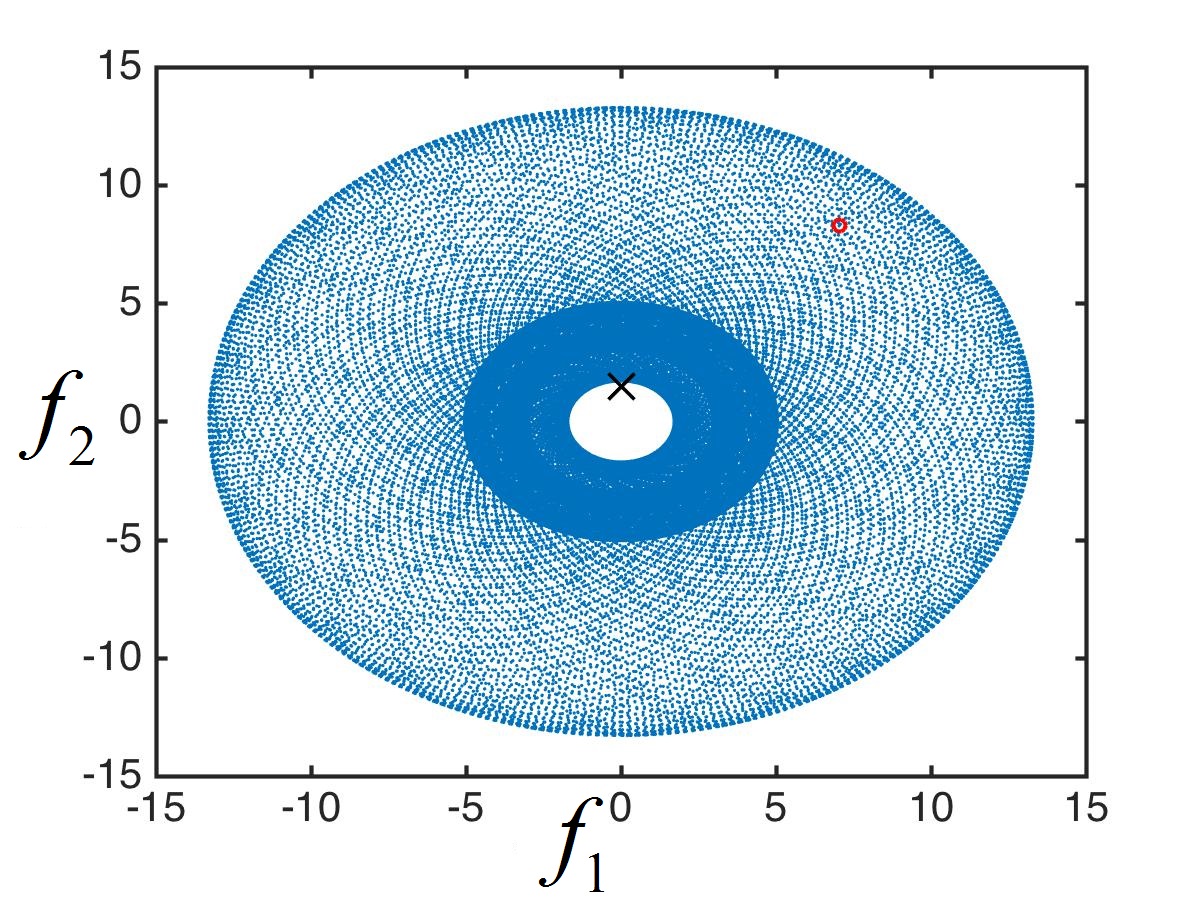}
\includegraphics[width=.45\textwidth]{\Path 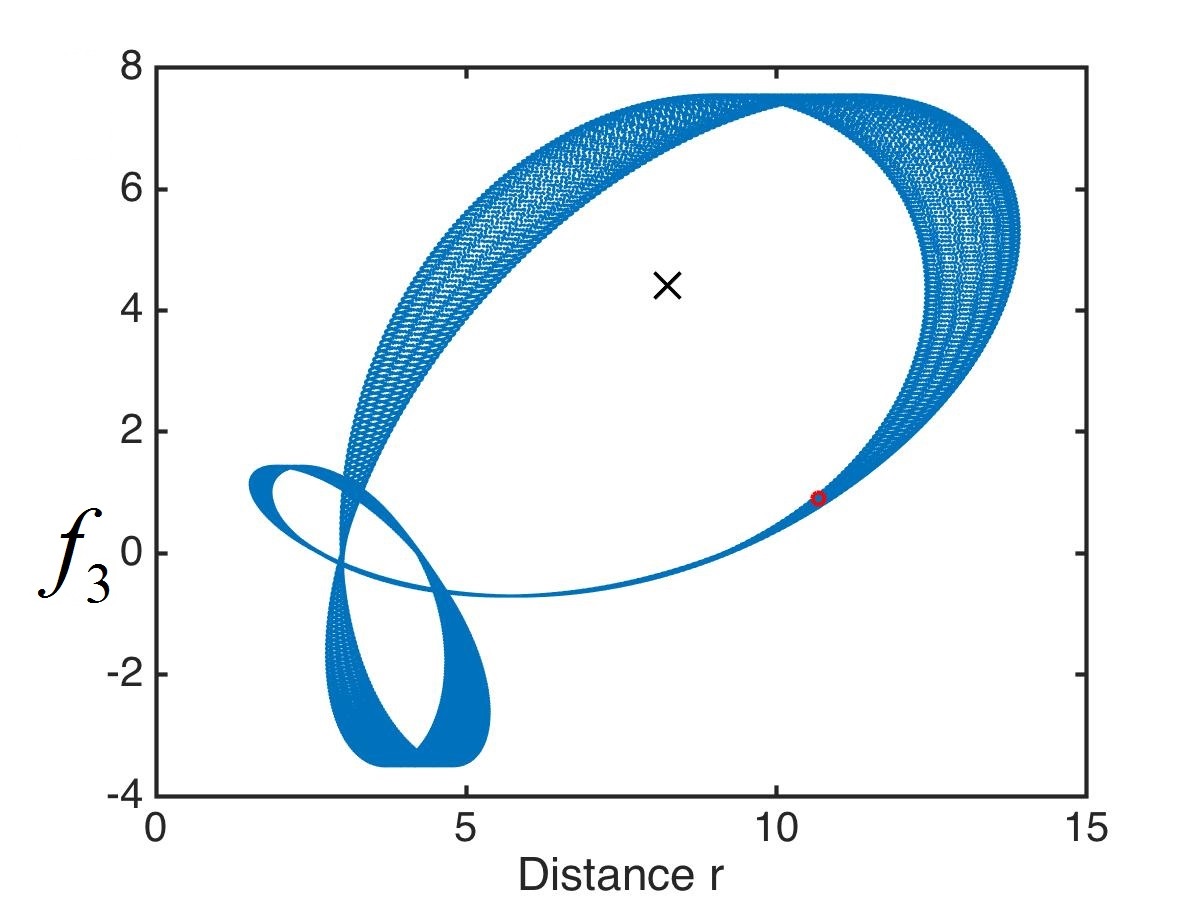}\\
\includegraphics[width=.45\textwidth]{\Path 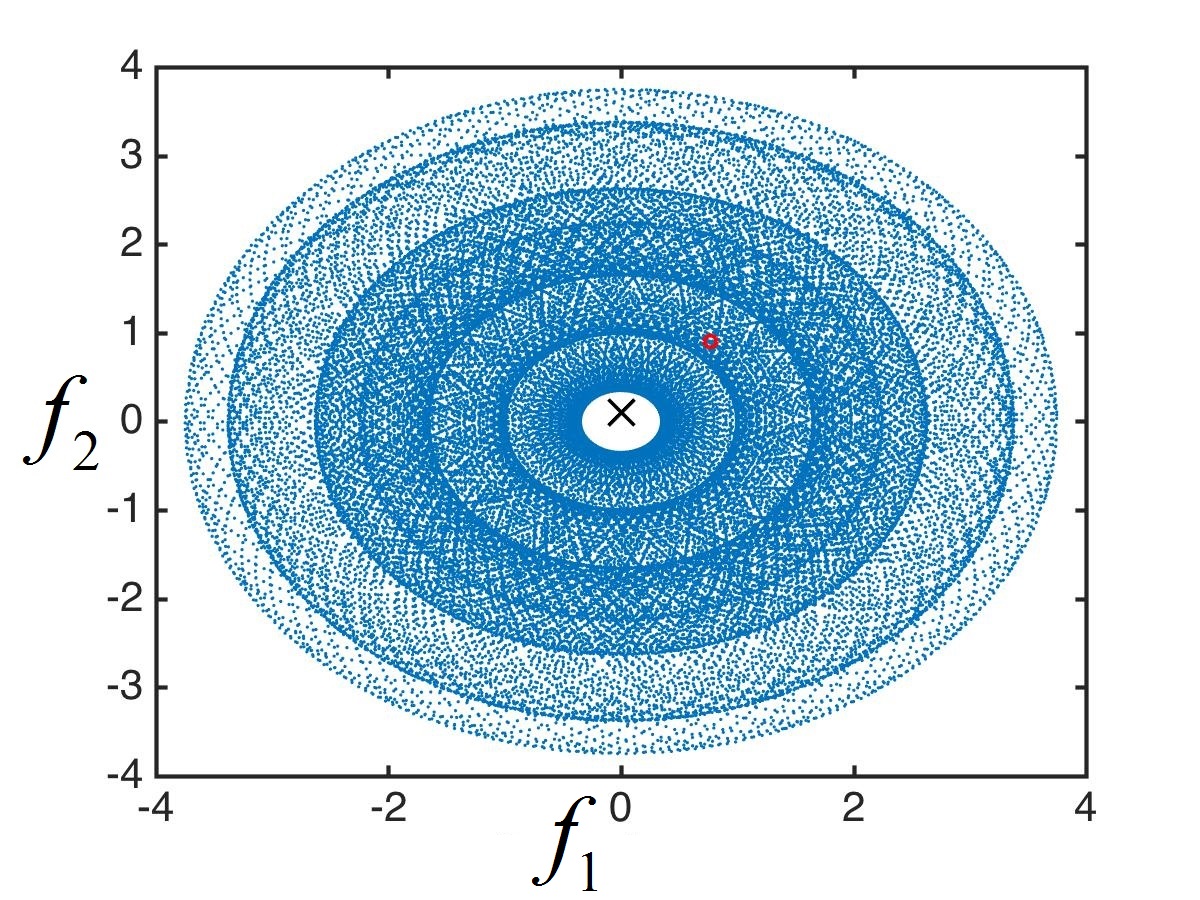}
\includegraphics[width=.45\textwidth]{\Path 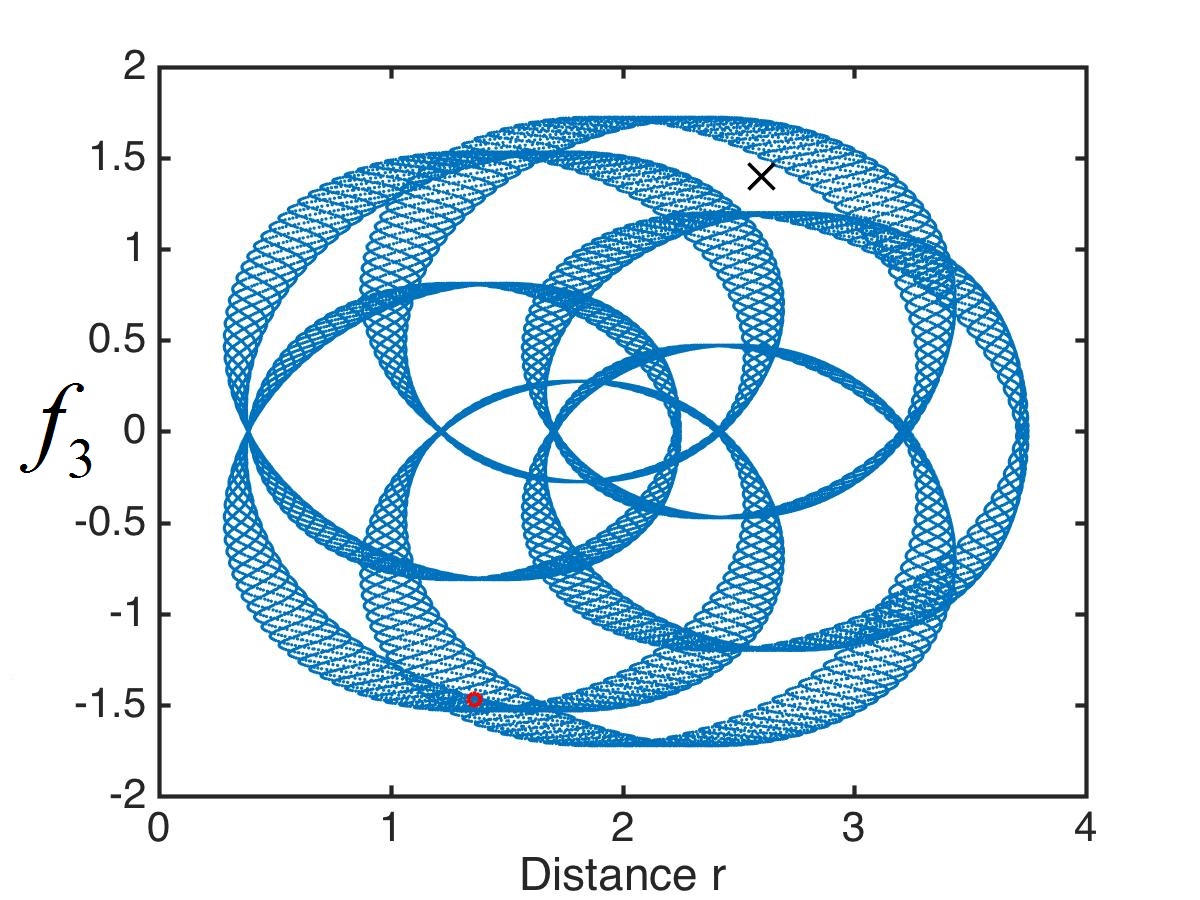}
\caption{{\bf Projections of the fish torus and the flower torus.} The coordinates used to find angle 1 (left) and angle 2 (right) for the fish torus (top) and the flower torus (bottom). The red circle shows the initial condition. The $\times$ shows the point with which the angle is measured. Note that for the 
the fish torus, the point from which the angle is measured is very close to the edge torus image. For angle 2, points are projected onto a tilted plane that makes angle $0.05 \pi$ with the horizontal. See Section \ref{subsec:2tori} for a full description of these projections.}
\label{fig:flw_fish_coordinates}
\end{figure}

\begin{figure}
\includegraphics[width=.45\textwidth]{\Path 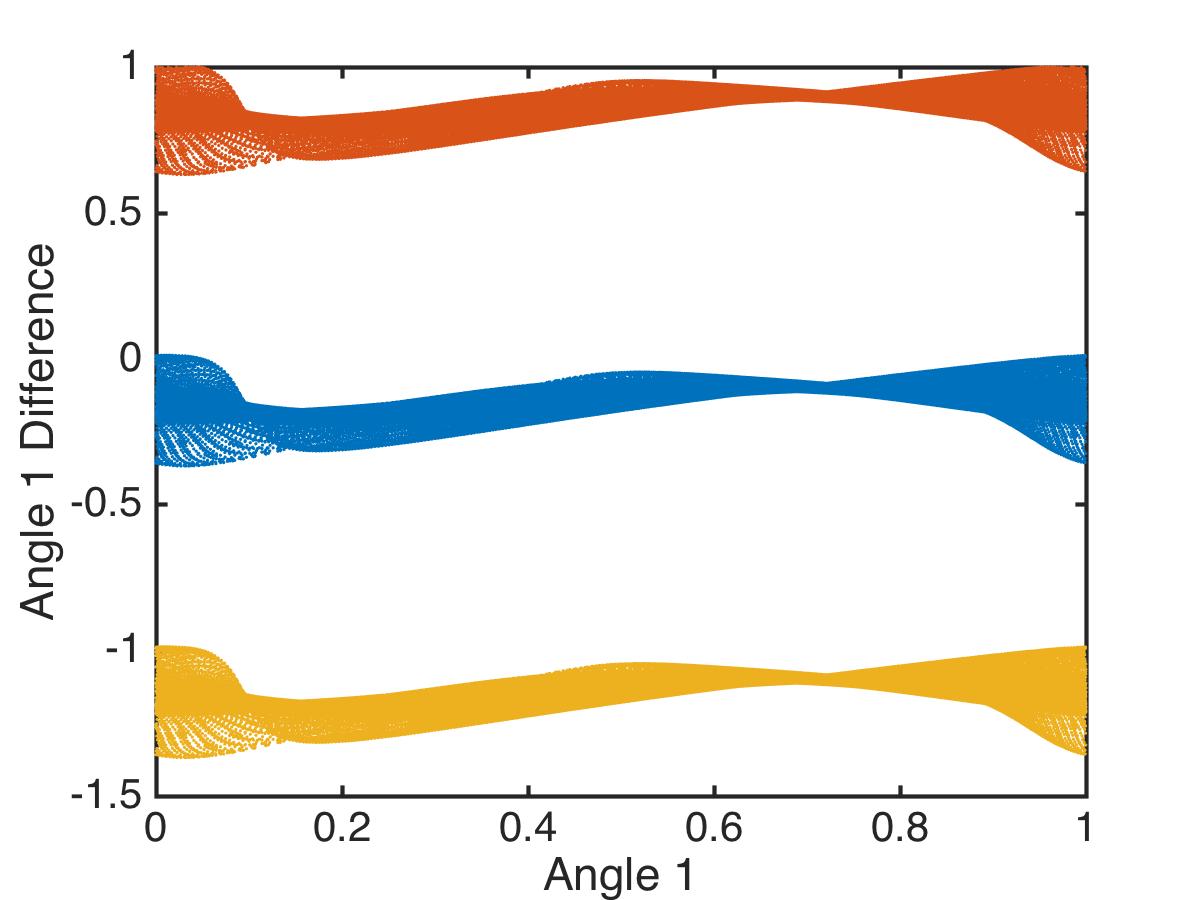}
\includegraphics[width=.45\textwidth]{\Path 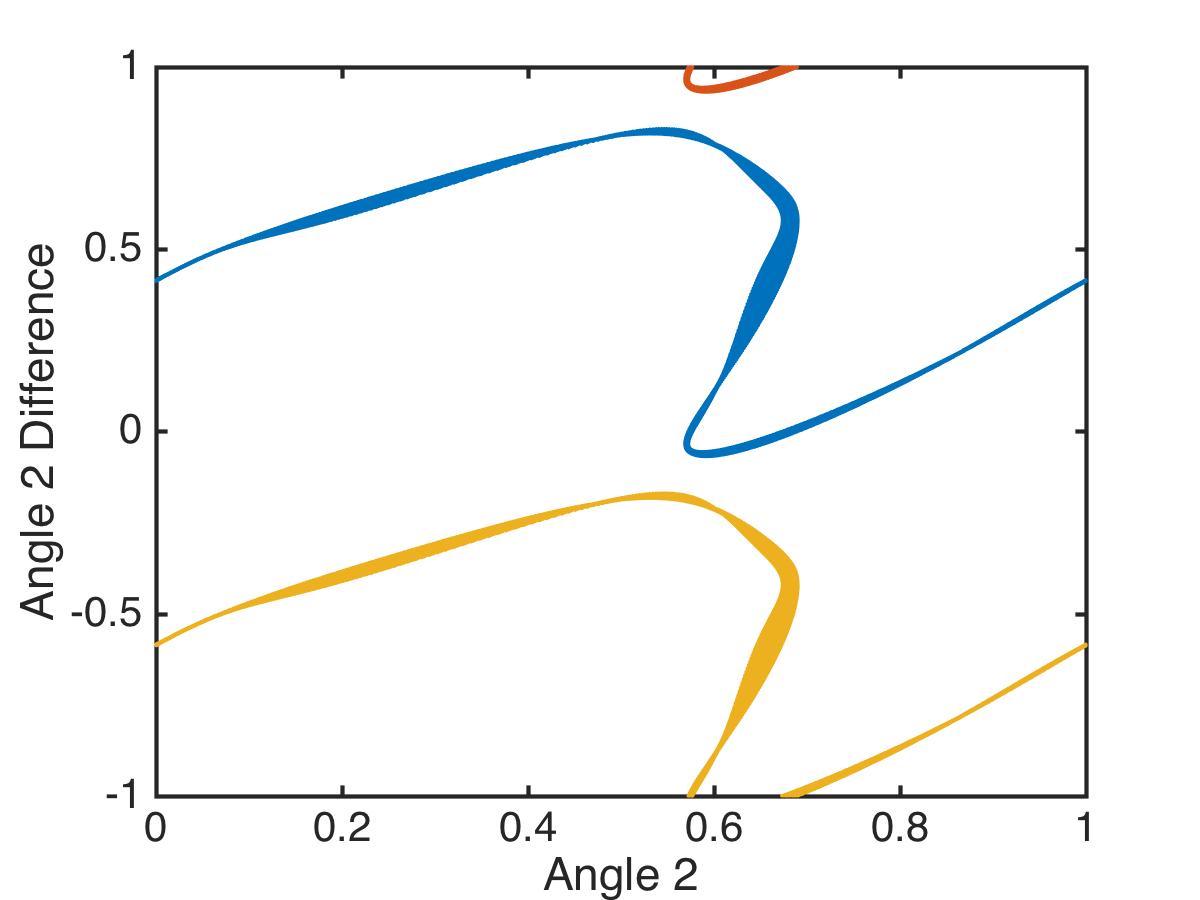}\\
\includegraphics[width=.45\textwidth]{\Path 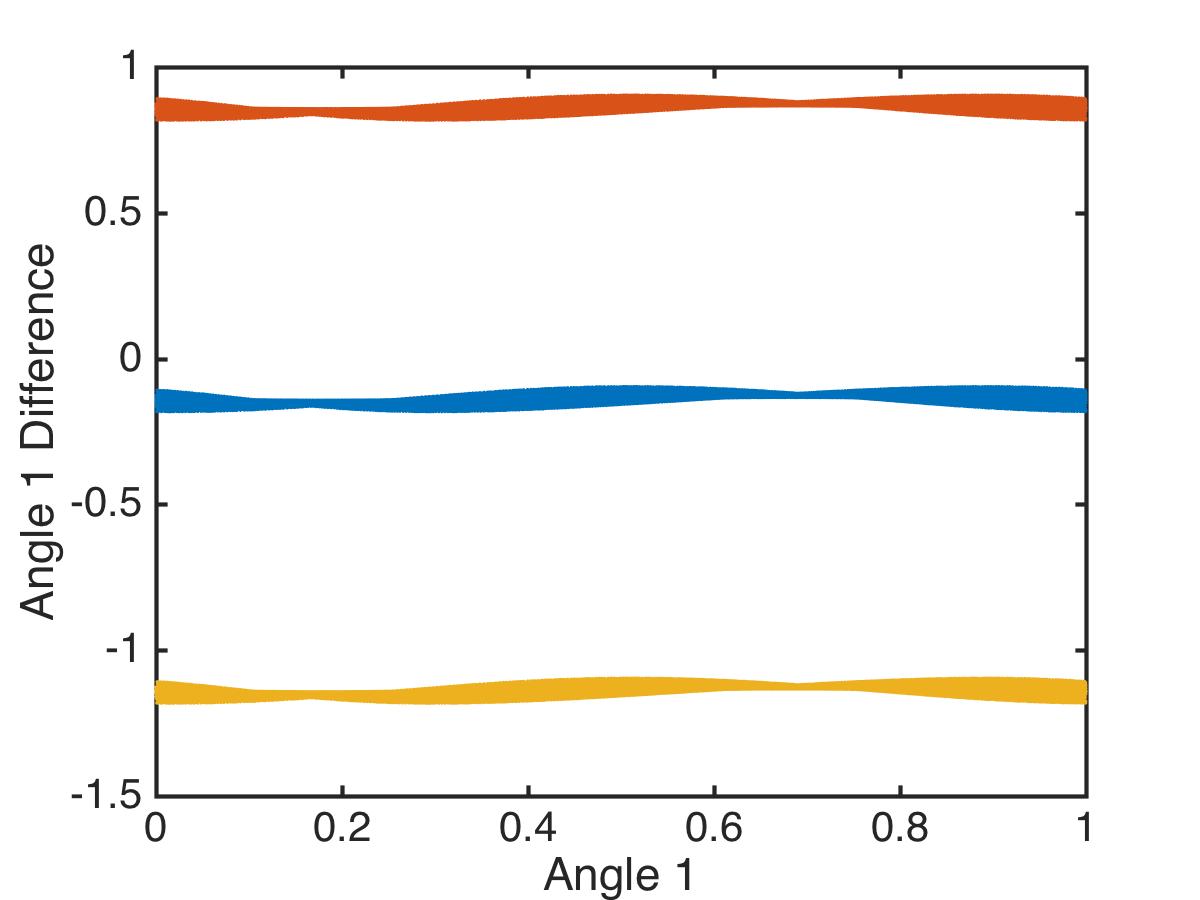}
\includegraphics[width=.45\textwidth]{\Path 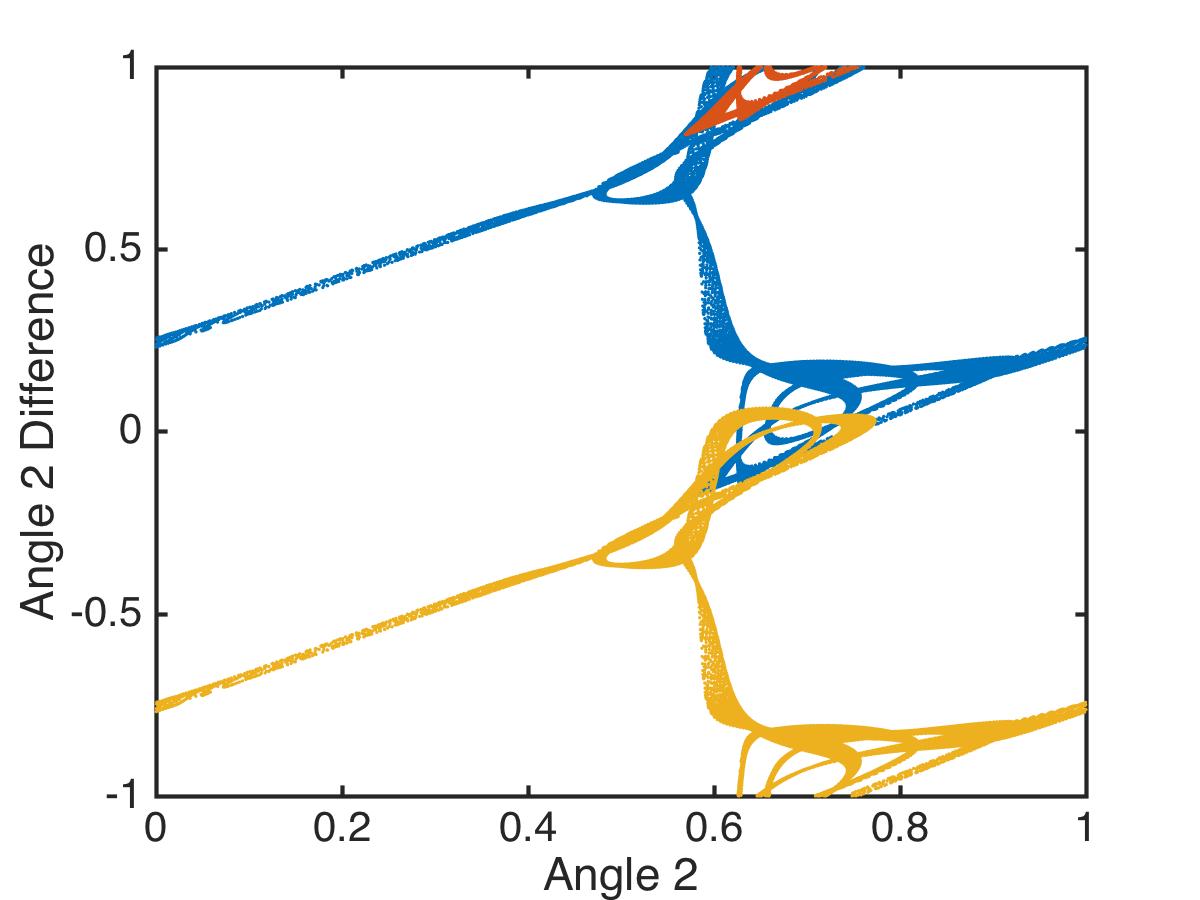}
\caption{{\bf Angle differences for the fish torus and flower torus.} 
Each panel shows three possible angle differences, each differing by an integer,
for the same projections as were depicted 
in Fig. \ref{fig:flw_fish_coordinates}.
The angle versus angle difference for angle 1 (left) and angle 2 (right) for the fish torus (top) and flower torus (bottom). 
In the final panel, the picture cannot be separated into separate components.
}
\label{fig:flw_fish_unlift}
\end{figure}

\begin{figure}
\includegraphics[width=.45\textwidth]{\Path 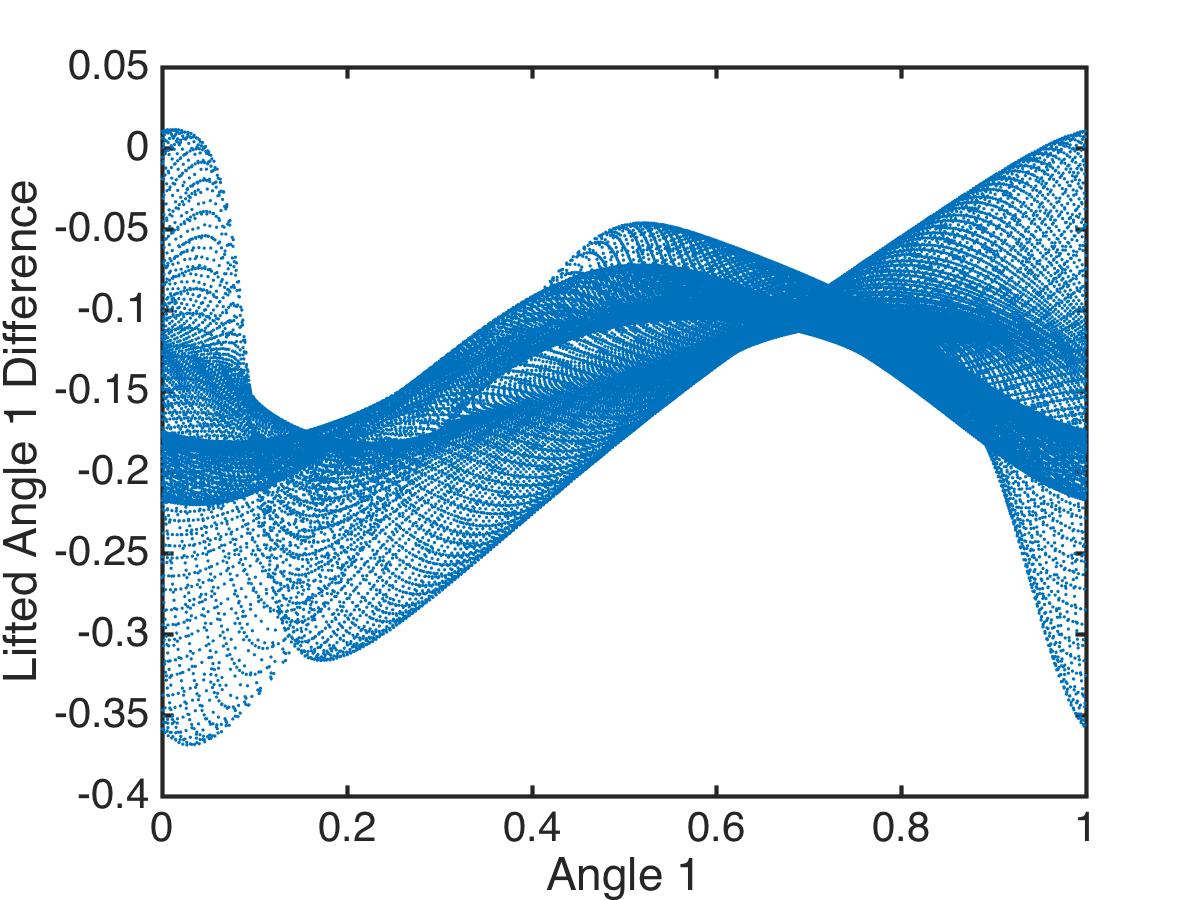}
\includegraphics[width=.45\textwidth]{\Path 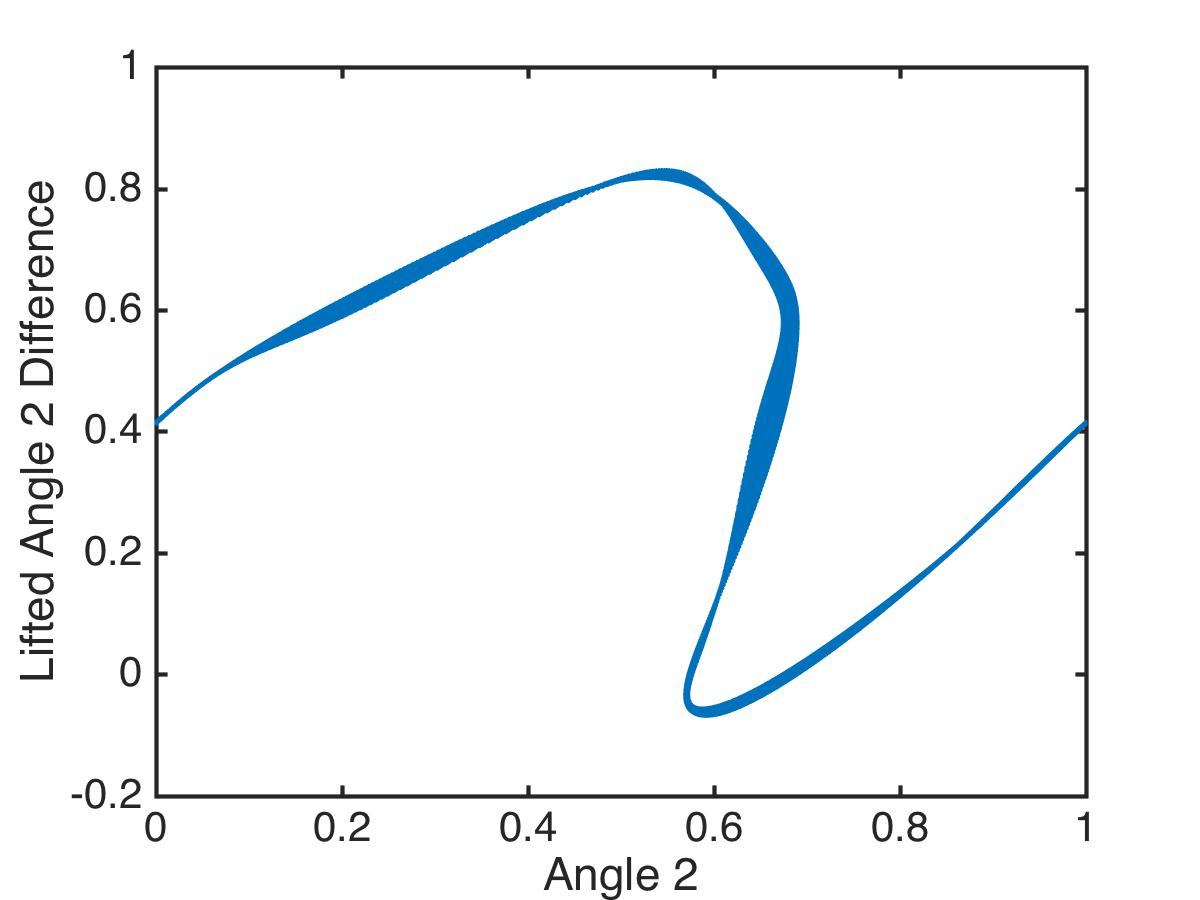}\\
\includegraphics[width=.45\textwidth]{\Path 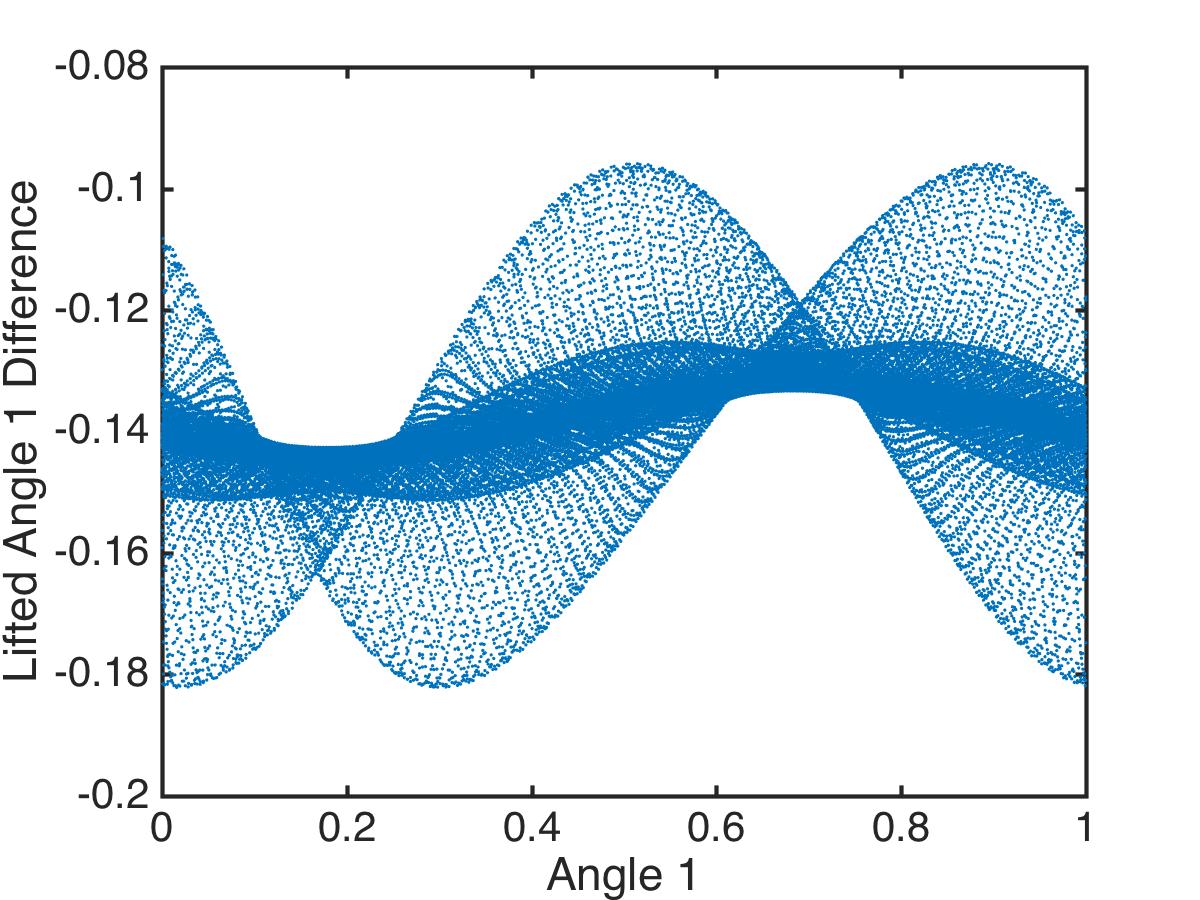}
\includegraphics[width=.45\textwidth]{\Path 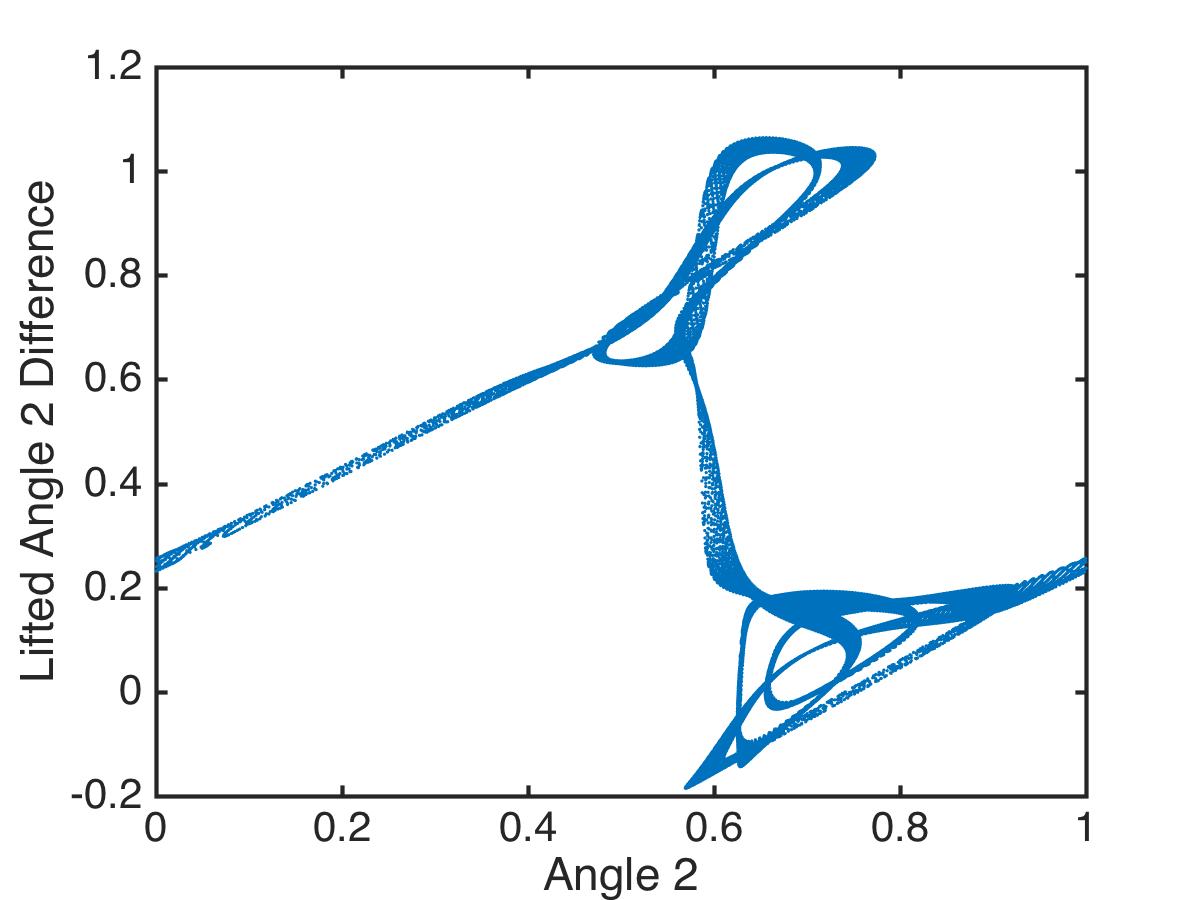}
\caption{{\bf Lifts of the angle difference for the fish torus and flower torus.} 
Here one of the possible lifts has been selected from each panel in Fig. \ref{fig:flw_fish_unlift}.
Each panel shows the angle versus angle difference lift for fish torus angle 1 (top left) and angle 2 (top right) 
and the flower torus angle 1 (bottom left) and angle 2 (bottom right), using the projections depicted 
in Fig. \ref{fig:flw_fish_coordinates}.}
\label{fig:flw_fish_lift}
\end{figure}

\subsection{Two examples in a higher dimension: fish and flower tori $\Torus$. \label{subsec:2tori}}
We use the fish and flower maps from the previous section in 
order to create 2-dimensional torus maps. We will explore the problem of computing rotation rates for these examples where we know the rotation rates for the quasiperiodic maps. 
Let $\rho := (\sqrt{5}-1)/2$ and $\phi := \sqrt{3}/2$, and define
\begin{eqnarray}
(\theta_n,y_n) := (n \rho~\bmod1, n \phi~\bmod1)\in \Torus
\end{eqnarray}
Let $\gamma$ be either the fish or the flower map defined in the previous section. 
Define the torus-version $f_T$ of the $\gamma$ map(s) as follows. Let
$Re(\cdot)$ and $Im(\cdot)$ denote the real and imaginary components of a complex number, and let $f_{T}:\Torus\to\R^3$. 
Write $f_T(\theta_n,y_n)= (f_1,f_2,f_3)(\theta_n,y_n)$, where
\begin{eqnarray}
f_1(\theta_n,y_n) &=& Re(\gamma(\theta_n) +2) \cos(2 \pi y_n)\\
f_2(\theta_n,y_n) &=& Re(\gamma(\theta_n) +2) \sin(2 \pi y_n)\\
f_3(\theta_n,y_n) &=& Im(\gamma(\theta_n)).
\end{eqnarray}
The ``+2'' is just for convenience so that the torus image can wrap around the 
origin rather than having to wrap it around some other point. 
For each $\gamma$, the map $f_T$ takes a quasiperiodic trajectory into $\R^3$. 

{\bf Two projections of a torus for two rotation rates.}
Figure~\ref{fig:flw_fish_coordinates} shows two independent projections of $f_T$ to $\R^2$. 
For the first rotation rate, we project $f_T$ to $(f_1,f_2)$ in the plane. Then we measure the angle $\phi$ from 
a reference point $P$ which is not in the image of the torus. In particular, $P=(0,1.5)$ for the fish torus, and $ (0,0.1)$ for the flower torus.
For both maps, this projection gives a rotation rate of $\phi/2 \pi$ (the denominator $2 \pi$ comes from the fact that we are measuring angles in $[0,1]$).

For a second rotation rate, let $R_\alpha$ be the rotation matrix that tilts by angle $\alpha = 0.05 \pi$ in the $f_2-f_3$ plane. Namely \[R_\alpha = \left(
\begin{array}{rrr}
1 & 0 & 0 \\
0 & \cos \alpha & -\sin \alpha\\
0 & \sin \alpha & \cos \alpha
\end{array}
\right)
\]
Set 
\[h = R_{0.05 \pi} f.\]
Define $r = \sqrt{ h_1^2 + h_2^2}$. Then our projection is to the value $(r,f_3)$. We measure the angle of this projection relative to the point 
$(8.25,4.4)$ for the fish torus, and $(2.6,1.4)$ for the flower torus. For both maps, this projection gives rotation rates of 
$1-\phi/(2 \pi)$
and $1-\rho$.
Why the tilt by $0.05 \pi$ rather than use value of $r$ with respect to the original coordinates? 
Because without the tilt (i.e. $\alpha = 0$), the projection would be a curve rather than a thick strip, which would not give a true test of our Embedding continuation method in two dimensions. 

In both cases, we get a map whose image has at least one hole (in which the winding number $=\pm 1$), and we can measure angles $\phi$ and angle differences
$\Delta$ compared to a point inside one of the holes, as long as the torus has a winding number $|W(p)|=1$
with respect to points in this hole. Thus just as for the one-dimensional case, we 
compute the lift, and then compute the rotation rate for these two different projections. 

Figures~\ref{fig:flw_fish_unlift} and \ref{fig:flw_fish_lift} show the original values of the angle difference, and 
the computed lift, respectively. 
Note 
that fish torus lift is easy to compute while the flower torus requires and embedding. 

As mentioned in Section \ref{sec:introduction}, rather than using Birkhoff Averages, we achieve more rapid convergence using our Weighted Birkhoff Average, denoted $\Q^{[p]}_N$ in Eq.\ref{eqn:QN}. 
Define the $\rho$ approximation $\rho_N := \Q^{[1]}_N(\hat\Delta_n)$, when $p=1$. Fluctuations in $\rho_N$ fall below $10^{-30}$ for $N>20,000$. 
Since we know the actual rotation rate, we can report that the error $|\rho-\rho_N|$ is then below $10^{-30}$.

\subsection{The circular planar restricted three-body problem ({\bf CR3BP})} \label{sec:3BP}
CR3BP is an idealized model of the motion of a planet, a moon, and an asteroid governed by Newtonian mechanics
Poincar{\'e}~\cite{ThreeBody2,ThreeBody1} introduced his method of return maps using this model. 
In particular, we consider a circular planar three-body problem consisting of two massive bodies (``planet" and a large ``moon") moving in circles about their center of mass and a third body (``asteroid") whose mass is infinitesimal, having no effect on the dynamics of the other two.

This model can also (simplistically) represent the Sun-Earth-Moon system discussed in the introduction though the parameter $\mu$ has to be changed, and the Moon is the body that is assumed to have negligible mass. All three travel in a plane.

We assume that the moon has mass $\mu$ and the planet mass is $1-\mu$ where $\mu = 0.1$, and
writing equations in rotating coordinates around the center of mass. Thus the planet remains fixed at $(q_1,p_1) = (-0.1,0)$, and the moon is fixed at $(q_2,p_2)=(0.9,0)$. In these coordinates, the satellite's location and velocity are given by the {\em generalized position vector} $(q_1,q_2)$ and {\em generalized velocity vector} $(p_1,p_2)$.

Define the distance of the asteroid from the moon and planet are 
\[ d_{moon}^{\ 2} = (q_1-1+\mu)^2+q_2^2\] \[d_{planet}^{\ 2} = (q_1+\mu)^2+q_2^2.\]
The following function $H$ is a Hamiltonian (see \cite{ThreeBody3} p.59 Eqs. 63-66) for this system
\begin{equation}\label{eqn:Hamiltonian}
H=\frac{1}{2}(p_1^2+p_2^2)
+p_1q_2-p_2q_1
-\frac{1-\mu}{d_{planet}} - \frac{\mu}{d_{moon}}, 
\end{equation}
where $p_1=\dot{q_1}-q_2$ and $p_2=\dot{q_2}+q_1$.
We get the equations of motion from 
\begin{equation*}
{\displaystyle
\begin{array}{rcl}
\cfrac{dq_i}{dt} &=& H_{p_i}, \\
\cfrac{dp_i}{dt} &=& -H_{q_i}. \\
\end{array}
}
\end{equation*}
That is, the equations of motion are as follows: 
\begin{equation*}
{\displaystyle
\begin{array}{rcl}
\cfrac{dq_1}{dt} &=& p_1+q_2, \\
\cfrac{dq_2}{dt} &=& p_2-q_1, \\
\cfrac{dp_1}{dt} &=& p_2- \mu\cfrac{q_1-1+\mu} {d_{moon}^{\ 3}} -(1-\mu)\cfrac{q_1+ \mu}{ d_{planet}^{\ 3}},\\
\cfrac{dp_2}{dt} &=& -p_1
-\mu\cfrac{ q_2} {d_{moon}^{\ 3}}
-(1-\mu)\cfrac{q_2}{ d_{planet}^{\ 3}},
\end{array}
}
\end{equation*}
\begin{figure}
\begin{center}
\includegraphics[width=.31\textwidth]{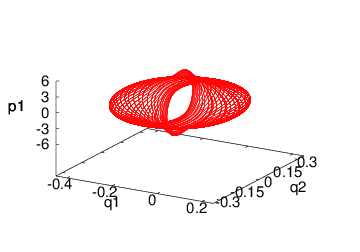}
\includegraphics[width=.36\textwidth]{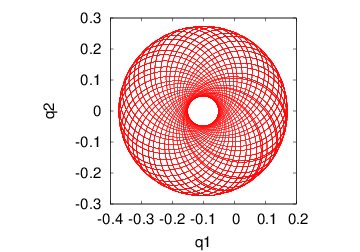}
\caption{Two views of a two-dimensional quasiperiodic trajectory for the restricted three-body problem described in 
Section~\ref{sec:3BP}.}
\label{fig:proj-b1}
\end{center}
\end{figure}

\begin{center}
\begin{figure}
\includegraphics[width=.49\textwidth]{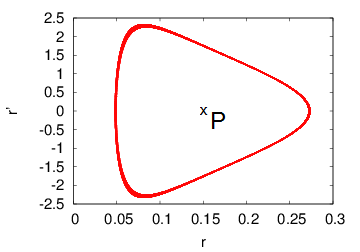}
\includegraphics[width=.49\textwidth]{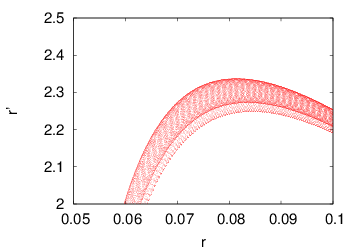}
\caption{\textbf{Plots of the circular planar restricted three-body problem in $r-r'$ coordinates.}
As described in the text, we define $r=\sqrt{(q_1+0.1)^2+q^2_2}$ and $r' = dr/dt$. 
This figure shows $r$ versus $r'$ for a single trajectory. The right figure is the enlargement of the left.
One of the two rotation rates $\rho^{*}_\phi$ is calculated by measuring from $(r,r^{\prime})=(0.15,0)$ in these coordinates.}
\label{fig:proj-r-rprime}
\end{figure}
\end{center}
\begin{center}
\begin{figure}
\includegraphics[width=0.49\textwidth]{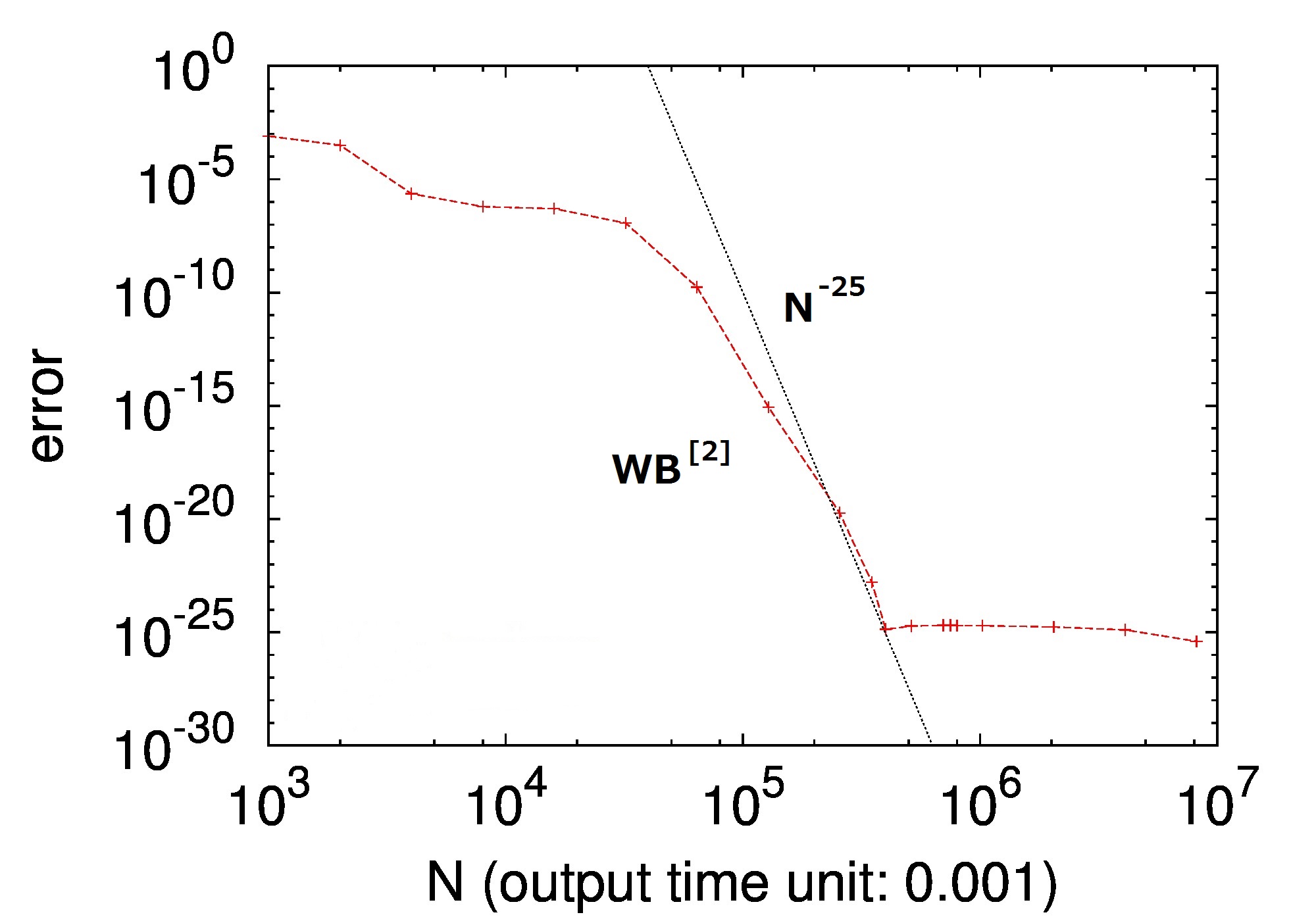}
\includegraphics[width=0.49\textwidth]{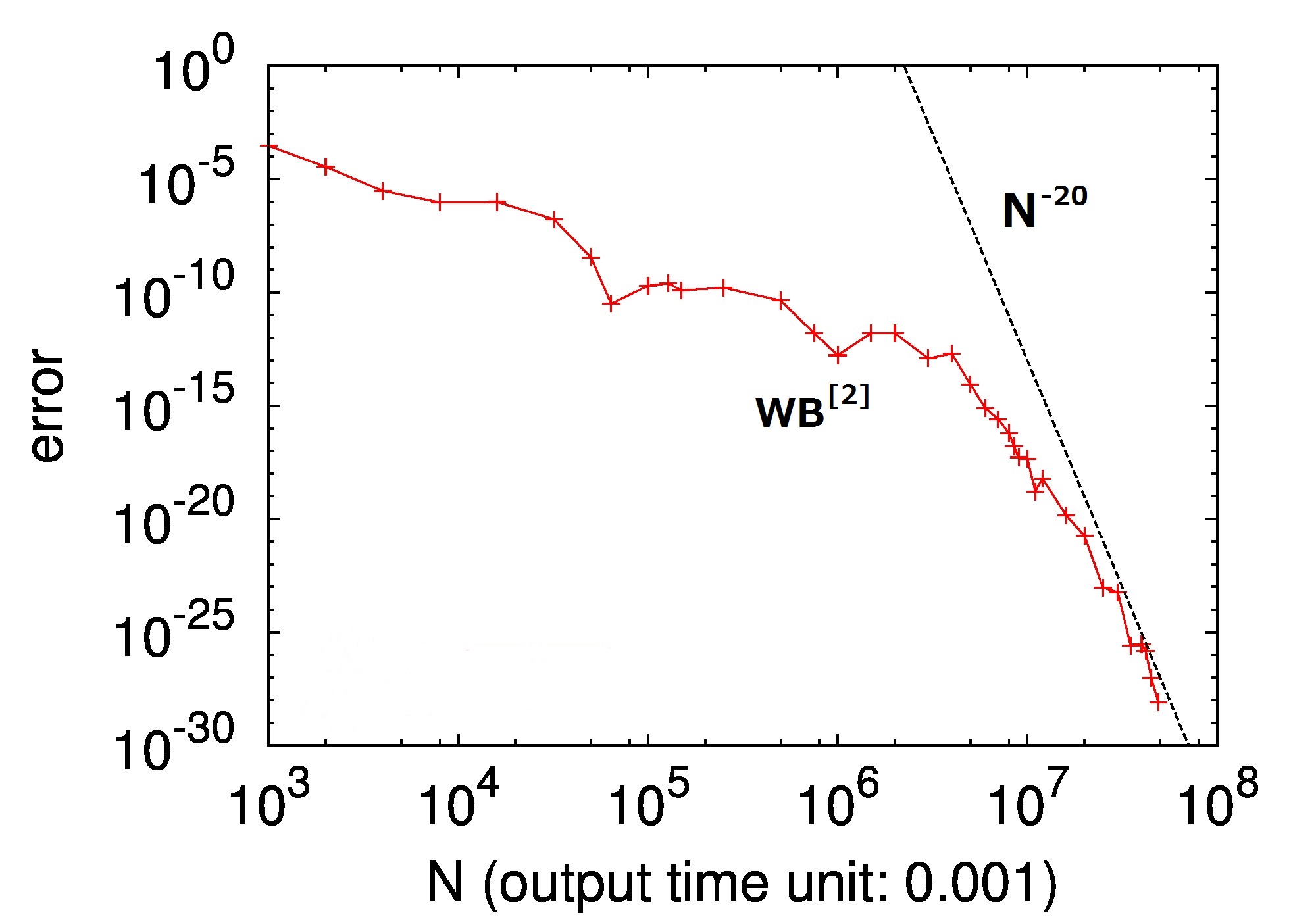}
\caption{\textbf{Convergence to the rotation rates for the CR3BP.} 
For these two figures, we used differential equation time step
$dt=0.00002$ and we compute the change in angle after 50 such steps, that is, in time ``output time'' $Dt=0.001$. 
We show the convergence rates to the estimated rate of $0.001\times \rho^{*}_{\theta}$ (left) and 
of $0.001 \times \rho^{*}_{\phi}$ (right).
For both cases rotation rates are calculated using the Weighted Birkhoff averaging method $\Q^{[2]}_{N}$ in Eq.\ref{eqn:QN} and show fast convergence. 
}
\label{fig:convergence-to-rotation-number}
\end{figure}
\end{center}
We measure angles as a fraction of a full rotation and not in terms of radians.
The asteroid's orbit in rotating coordinates is shown in 
Fig. \ref{fig:proj-b1}. Here time is continuous so we can measure the total angle through which a trajectory travels, retaining the integer part. 
The first rotation rate $\rho^{*}_{\theta}$ of the asteroid's orbit is its average rate of rotation about the planet, that is, the average rate of change of the angle $\theta$ measured from $(q_1,q_2)=(-0.1,0)$. We compute that $\rho^{*}_{\theta}
=-2.497823504839344460408394$ rev/sec, 
that is, about -2.5 $\theta$-revolutions per unit time 
Fig. \ref{fig:convergence-to-rotation-number} (left) shows the error in convergence to the value
0.001$\rho^{*}_{\theta}$ 
Note that the rotation rate in the fixed coordinate frame is $\rho^{*}_{\theta} + 1/(2 \pi).$

The second rotation rate $\rho^{*}_{\phi}$ measures the oscillation in
the distance $r$ from the planet. In particular, we project to
the $(r,r^{\prime})$ plane, where $r^{\prime} := dr/dt$. 
That is, define
$r=\sqrt{(q_1+0.1)^2+q^2_2}$ and
$r^{\prime}=\frac{dr}{dt}=((q_1+0.1)\frac{dq_1}{dt}+q_2
\frac{dq_2}{dt})/r$, as shown in Fig. \ref{fig:proj-r-rprime}. The angle $\phi$ is
measured from $(r,r^{\prime})=(0.15,0)$. The fast convergence to the
value 0.001$\rho^{*}_{\phi}$ by the Weighted Birkhoff Average $\Q^{[2]}_{N}$ in \ref{eqn:QN} is seen in
Fig. \ref{fig:convergence-to-rotation-number} (right), where
$\rho^{*}_{\phi}=-2.3380583953388194764236520190142509$ rev/sec. The period
of time between perigees is the reciprocal, or about 0.43 time units.

We used the 8th-order Runge-Kutta method in Butcher~\cite{butcher08} to compute trajectories of CR3BP with time steps of $h=2\times 10^{-5}$.

\textbf{The meaning of rotation rates for the CR3BP.} 
In \cite{DSSY}, we investigated the same asteroid orbit of the CR3BP as is studied here, but instead of 
the continuous-time trajectory that lies on a two-dimensional torus as presented above, 
there we used a Poincar\'e map. The coordinates of the asteroid were recorded each time 
the asteroid crossed the line $q_2=0$ with $dq_2/dt>0$. In the cases we study, the map trajectory is a quasiperiodic trajectory on a closed curve. 
Hence there is only one rotation rate, a much simpler situation. 
Choosing a point inside the closed curve, we computed a rotation rate, namely the average angular rotation {\it per iteration of the Poincar\'e map}. The rotation rate $\rho^{*}_{P}$ per Poincar\'e map on the Poincar\'e surface $q_2=0$ (or equivalently, $\theta=0$) around $(q_1,p_1)=(-0.25,0)$ was computed as 0.0639617287574530971640777244014426955. 
We felt that the issues of rotation rates could be clarified if we computed the trajectory as a continuous orbit as we do here. The two rotation rates computed here
$\rho^{*}_{\phi}$ and $\rho^{*}_{\theta}$ and our previous result $\rho^{*}_{P}$ bear the following relation to our previous results:
\[\rho^{*}_{P}=\left( \pm \cfrac{\rho^{*}_{\phi}}{\rho^{*}_{\theta}} \right) \bmod 1. \]
\[
\rho^{*}_{P} = 0.06396\cdots = 1 - \cfrac{2.338\cdots}{2.497\cdots} \pm 10^{-25}
\]
See the caption of Fig. \ref{fig:convergence-to-rotation-number}. We solved the differential equation using an $8^{th}$-order Runge-Kutta method using quadruple precision. Both approaches are based on rotating coordinates, but there is another approach.

{\bf The orbit as a slowly rotating ellipse.} The asteroid rotates about the planet at a rate of $\rho^{*}_\theta$ revolutions per unit time when viewed in the rotating coordinate in which the moon and planet are fixed. The sidereal rotation rate (as viewed in the coordinates of the fixed stars) is $\rho^{*}_\theta\ + 1/(2 \pi)$.
We can think of the orbit as an approximate ellipse whose major axis rotates and even changes eccentricity (being more eccentric when the asteroid apogee is aligned with the planet moon axis). 

Without the moon the asteroid orbit would be perfectly elliptical with its major axis fixed in position, but the moon causes the ellipse to rotate slowly.
The angle $\phi(t)$ tells where the asteroid is on its roughly elliptical orbit; Fig. \ref{fig:proj-r-rprime} shows that the apogee occurs when when the distance from the planet $r \sim 0.27$ and the perigee when $r\sim 0.05$, with some variation. The time between successive perigees averages $1/\rho^{*}_\phi$. Note that the difference in these rates satisfies
$$ \rho^{*}_\phi \ - [\rho^{*}_\theta \ + 1/(2 \pi)] \sim 0.000610166 \sim 1/1638.9 .$$ Hence relative to the fixed stars, that is, in non-rotating coordinates, the asteroid's ellipse's major axis precesses slowly. Its apogee point returns to its original position (in non-rotating coordinates) after the asteroid passes through its apogee approximately $1639$ times. 


\section{Discussion and conclusions}\label{discussion}

What does it mean to ask for one or more rotation rates of the $d$-dimensional quasiperiodic map Eq. \Xx? One might expect that one should find $\rho$ or rather its coordinates. As we explain below and in Section~\ref{sec:rotvecset},
this is an ill-posed problem (especially for $d>1$. The Babylonians computed rotation rates for projections of the Moon's trajectory onto the globe of fixed stars (as we have discussed in the Introduction). So we refer to their approach as the ``Babylonian Problem'': computing rotation rates for a projection of a quasiperiodic process.

We have developed our Embedding continuation method for calculating the rotation rate for ``almost every'' 
Babylonian Problem, that is, for
smooth projection $\psi$ from of a quasiperiodic dynamical system on $\torus$. ``Almost every'' is in the sense of prevalence - and in practice there will be difficult cases especially since the number $N$ of interates needed increases as $d$ increases. Our Weighted Birkhoff Method of computing rotation numbers significantly shortens the computation time for computing rotation numbers, making our approach effective in practice. See the Introduction. 

A key motivating difference between $d=1$ and $d>1$ is that in the higher-dimensional case, for the rotation vector $\rho\in\torus$ there are infinitely many ways of choosing coordinates on $\torus$ for the map Eq. \Xx. 
In Section~\ref{sec:rotvecset} we show that 
the set of resulting coordinate representations $(\rho_1,\cdots,\rho_d)$ of $\rho$ are dense in $\torus$. 
Every point $r$ in $\torus$ is arbitrarily close to such representations. Hence instead of trying to find the coordinates of $\rho$, we have learned from the Babylonians, and we phrase our goals in terms of finding a rotation number $\rho_\psi$ (usually, $\rho_\phi$ or $\rho_\gamma$) for some projection from $\torus$ into a one or two-dimensional space.

Even for $d=1$, there is some uncertainty for obtaining $\rho$  depending on the choice of orientation on $S^1$. 
We can obtain either $\rho ~\bmod1$ or $1-\rho ~\bmod1$.

In Section~\ref{sec:3BP}, we apply our method to the quasiperiodic torus occurring for a 4-dimensional circular restricted 3-body problem, 
depicted in Fig.~\ref{fig:proj-b1}. In particular, we explain the relationship between the two rotation rates obtained 
from the original differential equation system and the rotation rate which was previously obtained from the Poincar\'e map. The fact that the rotation rate of an asteroid will be different depending on whether on uses rotating coordinates or sidereal coordinates (in which the distant stars are fixed) is an example of how the rotation rate can depend on the projection.

{\bf Notes on delay coordinate embedding theorems.}
H. Whitney \cite{Whitney} showed that a topologically generic smooth map $\Gamma$ from a $d$-dimensional smooth compact manifold $M$ into $\R^{D}$ where $2d+1\le D$ is a diffeomorphism on $M$; in particular the map $\Gamma:M\to F(M)$ is an embedding of $M$. 

Sauer et al  \cite{Embedology}
modified Takens' result in two ways. First, it replaced ``topologically generic'' by ``almost every’’ (in the sense of ``prevalence’’) 
in Theorems 2.1 and 2.3 in \cite{hunt}. See also \cite{ott}. For physical purposes ``almost every'' has significance while residual sets do not seem to.
In this paper, in Theorem \ref{thm:takens}, we have adapted the ``almost every'' approach.

For completeness, we mention the second way  \cite{Embedology} generalized Takens' approach, even though this second way is not used here, because the sets we deal with are manifolds. 
The second way is that 
 \cite{Embedology} 
 allblack
 replaced the assumption that $M$ is a manifold by assuming only that $M\subset \R^k$ for some $k$ is an invariant set of some map and that $M$ has box dimension $boxdim(M)$ and $\Gamma$ is a mapping of a neighborhood of $M$ into $\R^D$ where $D>2\cdot boxdim(M)$. 
The great majority of citations to Takens \cite{Takens} are for the case where $M$ is a chaotic attractor that is not a manifold so that Takens' Theorem does not apply.
Those papers actually use the results in \cite{Embedology}, not in Takens' \cite{Takens}. 
One unusual aspect of our current paper is that we actually only need the case that Takens proved. Here $M$ is a quasiperiodic torus so it is a manifold. 

The Takens Theorem also has assumptions that the set of periodic points $F:M\to M$ for some smooth map was in some sense small, in our case there are no periodic points so those assumptions are automatically satisfied. 
Hence we only state it in a special case needed here.

  \label{sec:discussion} We have demonstrated that in one dimension, a
rotation rate can be computed precisely with minimal ambiguity, but 
higher dimensional cases ($\torus$ with $d>1$) are more complicated. Projections into the plane
yield rotation rates, but there are infinitely many topologically distinct ways to
project a higher dimensional torus onto a circle, each of which yields a
different rotation rate. This makes it important for the investigator to
explain the meaning of any particular rotation rate. In fact, a rotation rate
is a rate specifying an average change per unit time,
where there can be considerable choice in the time units. To illustrate this point,
we more carefully consider the CR3BP example with a focus on what the rotation rates tell us about the trajectories of an asteroid.

{\bf Acknowledgments.} We would like to thank the referee for many helpful comments. 
YS was partially supported by JSPS KAKENHI grant 17K05360 and JST PRESTO grant JPMJPR16E5.
ES was partially supported by NSF grant DMS-1407087. 
JY was partially supported by National Research Initiative Competitive grants 2009-35205-05209 and 2008-04049 from the USDA.

\bibliographystyle{unsrt}
\bibliography{qr_bibliography,Weighted_calc_bibliography}
\end{document}